\documentclass{article}

\usepackage{natbib}%
\usepackage[figuresright]{rotating}

\usepackage{url}
\usepackage{algorithmic}
\usepackage{algorithm}

\usepackage{amsmath,amsthm,amssymb}

\theoremstyle{definition}
\newtheorem{example}{Example}

\theoremstyle{remark}

\usepackage{xcolor}
\definecolor{ForestGreen}{RGB}{34,139,34}
\usepackage{mathtools}

\usepackage{tikz}
\usepackage{caption}
\usepackage{subcaption}
\usepackage{multirow}

\usepackage{geometry}
\newcommand{\Day}{\operatorname{day}}
\newcommand{\Time}{\operatorname{time}}
\newcommand{\VLEN}{\operatorname{len^{\textsc{v}}}}
\newcommand{\TLEN}{\operatorname{len^{\textsc{t}}}}
\newcommand{\VCOST}{\operatorname{cost^{\textsc{v}}}}
\newcommand{\TCOST}{\operatorname{cost^{\textsc{t}}}}
\newcommand{\NN}{\mathbb{N}}

\newcommand{\vacio}{\circ}
\newcommand{\lleno}{\bullet}
\newcommand{\arcLTC}{\operatorname{arc_{\textsc{ltc}}}}
\newcommand{\arcLTR}{\operatorname{arc_{\textsc{ltr}}}}
\newcommand{\origen}{\operatorname{L}}
\newcommand{\destino}{\operatorname{L'}}
\newcommand{\inicio}{\operatorname{\mathcal{I}}}
\newcommand{\fin}{\operatorname{\mathcal{I}'}}
\newcommand{\capacity}{\operatorname{cap}}

\begin{document}

% Title
\title{Exact resolution of a simultaneous vehicle routing and crew scheduling problem in long-haul transport}

%Authors, affiliations address.
\author{Mauro Lucci\\
	\texttt{mlucci@fceia.unr.edu.ar}\\
	CONICET - Universidad Nacional de Rosario, Argentina.\\
	\and 
	Daniel Sever\'in\\
	\texttt{daniel@fceia.unr.edu.ar}\\
	Universidad Nacional de Rosario, Argentina.
	\and
	Paula Zabala\\
	\texttt{pzabala@dc.uba.ar}\\
	CONICET - Universidad de Buenos Aires, Argentina.}

\date{}

\maketitle

%Abstract
\begin{abstract}
This work focuses on exact methods for a Simultaneous Vehicle Routing and Crew Scheduling Problem in long-haul transport. Pickup-and-delivery requests with time windows must be fullfiled over a multi-day planning horizon. Unlike some classic approaches, the correspondence between trucks and drivers is not fixed and they can be exchanged in some locations and at any time. Drivers can also travel for free as truck passengers or take external taxis for an additional cost. The objective is to minimise the truck and taxi travel costs and the penalties for late deliveries. Routes for trucks and drivers are represented separately as directed paths in certain digraphs and then synchronised in time and space. Three compact Integer Linear Programming formulations are proposed and many families of valid inequalities are described. Extensive computational experiments are conducted on randomly generated instances. The formulations are experimentally compared and the effectiveness of the proposed valid inequalities as cutting planes in a branch-and-cut algorithm is evaluated.
\end{abstract}

%Keywords, etc.
\textbf{Keywords}: Integer Programming; Valid inequalities; Long-haul transport; Vehicle routing; Crew scheduling.

% Heading 1
\section{Introduction}\label{sec:Introduction}

Vehicle Routing Problems (VRPs) are among the most popular in transportation and logistics.
Applications go from supplying fuel to gas stations, collecting fresh produce from farms, long-haul and last-mile transport, home health care, and dial-a-ride services, to many others \citep[see][]{Toth2014}.
The opportunity to minimise the high costs they entail in the supply chain and their notorious hardness from the point of view of combinatorial optimization explains the broad interest they have aroused.
New requirements that constantly emerge in real-world applications often limit or prohibit the practical use of commercial software for general-purpose VRPs, motivating the development of specific approaches.

In recent decades, there has been growing interest in considering hours-of-service regulations for drivers during vehicle routing.
They are generally regulated by labour laws, collective labour agreements, and internal company policies.   
In this work, we focus on long-haul transport, whose requirements differ greatly from other transport sectors such as airways, railways, and urban mass transit, since trips are not timetabled and can be interrupted for drivers to take breaks.
Considering that the truck travel times can extend to several days, taking driver rest periods into account is essential to building feasible routes.
There is extensive literature in this regard, e.g. for United States law see \citet{Rancourt2013, Goel2014, Koc2018}, for European Union law see \citet{Goel2009, Kok2010, PrescottGagnon2010, Goel2014}, and for exact methods see \citet{Goel2017, Tilk2020}. 
For crew scheduling in other transport sectors see 
\citet{DEVECI2018, HEIL2020, IBARRAROJAS2015}.
It is also worth mentioning the work of \citet{goel2020team}, who compare the advantages of team vs. single driving in European road freight transport.
For that social legislation, different laws apply to team driving, e.g. the total daily driving time is doubled and a driver can take a break while the other is driving.
Generally, multiple objective functions are optimised simultaneously though the most common are the number of vehicles and drivers used, the total travel distance of the vehicles, and the total working time of the drivers.

All the above works for long-haul transport assume a fixed correspondence between trucks and drivers.
That is, the same driver travels the entire route of a truck, thus remaining unnecessarily idle during her/his rest periods.
There is room for improvement when drivers are allowed to change trucks over time; e.g. when a driver needs to take a break, another can relieve her/him and continue driving without delaying transport.
This motivates some recent works focused on Simultaneous Vehicle Routing and Crew Scheduling Problems (SVRCSPs), where such a fixed correspondence is relaxed.
For these variants, any route plan where drivers do not change trucks remains feasible, but more efficient use of resources can lead to new solutions to improve overall costs.
As a counterpart, the combinatorics grows since routes must also be built for drivers, which must be synchronised in time and space with the truck routes.
Besides, sinchronisation constraints introduce interdependence, meaning that changes to one route can affect the feasibility of the others \citep[see][]{Drexl2012}. 
In the brief literature on SVRCSPs, the most common resolution method consists of a two-stage sequential decomposition, e.g. see \citet{Drexl2013, Mendes2020decision, Lucci2022}.
In the first stage, restrictions on drivers are completely or partially ignored to build truck routes that serve all customers.
After that, drivers are assigned to segments of the previous routes in compliance with labour laws.

Although decomposition is effective in solving large instances in reasonable times, there is no guarantee of global optimality; e.g. crew scheduling could be too costly or even impossible for a given input (truck routes).
Furthermore, decomposition indirectly introduces a hierarchical order between the objectives, which makes it difficult to prioritise the objectives of the second stage over those considered in the first.
Exact methods that address the problem in a single stage can overcome these drawbacks, but they are usually more complex and more time-consuming.
To our knowledge, very few works in the literature deal with exact methods for SVRCSPs.
\citet{Lam2015} consider an application in humanitarian and military logistics, where crews can change vehicles in different locations and also travel as passengers, and the objective is to minimise a weighted sum that considers the number of vehicles and crews used and the total vehicle and crew travel distances.
The authors propose a constraint programming formulation that is solved using a LNS, a mixed integer program that is directly optimised by a general-purpose ILP solver, and a two-stage method.
From computational results on instances with up to 96 requests, they conclude that a single-stage method produces considerable beneﬁts over a sequential one, the constraint program scales better than the mixed integer program, and decoupling crews from vehicles is crucial to obtaining high-quality solutions.
\citet{DominguezMartin2018} consider an application inspired by local air traffic operations in the Canary Islands, which involves two depots. 
Vehicle routes start and end at different bases, but since crews always return to their starting base, drivers can switch vehicles in some exchange locations and travel as passengers. 
The authors present an ILP formulation to minimise the cost of driver routes, introduce some valid inequalities, and develop separation procedures.
They are then incorporated into a branch-and-cut algorithm to solve instances with up to 30 customer locations.

Both mathematical models proposed by \citet{Lam2015} and \citet{DominguezMartin2018} limit vehicle exchanges to occur only when customers are being served and also restrict drivers to single shifts.
Some recent works seem to reveal that greater savings can be achieved when drivers have more opportunities to switch vehicles.
Despite not involving a VRP, it is worth mentioning the work of \citet{AMMANN2023}, who address a driver routing and scheduling problem in long-distance bus networks, where drivers are allowed to exchange buses at arbitrary intermediate stops en route.
They define a MIP formulation based on a time-expanded multi-digraph and propose a destructive-bound-enhanced matheuristic.
A computational study with real-world data concludes that considering driver exchanges at intermediate stops in addition to regular stops allows additional cost savings of 43\%. % on average.

This work seeks to contribute to the study of exact single-stage methods for solving SVRCSPs, where the planning horizon extends over multiple days, trucks can have up to two drivers, and drivers can change trucks at any time in a set of locations.
A case study is followed that involves performing a set of pickup-and-delivery transport requests with multiple time windows over long distances, and the objective is to minimise a weighted sum that considers the travel costs and penalties for late deliveries.  
We propose three digraphs-based ILP formulations that model the movement of trucks and drivers over time, evaluate alternatives for some of the constraints, study some families of valid inequalities and their incorporation into branch-and-cut (B\&C) algorithms, and perform extensive computational experiments to compare the algorithms and evaluate their limits.
Some preliminaries were already presented in the conference paper \citet{Lucci2022}.
A coffee distribution company inspires the case study and was also addressed in \citet{Lucci2023} through a two-stage sequential decomposition.
In that previous work, hybrid metaheuristics were proposed to solve both stages, including a LNS with SA stopping criteria for the first stage and a GRASP$\times$ILS for the second.
Several generated instances were used to test the algorithms with up to 3000 requests distributed in 15 Argentine cities and a planning horizon of 1 to 4 weeks.

The rest of the paper is structured as follows. 
Section \ref{sec.problem.description} presents the study case.
Section \ref{sec.routes.representation} describes the representation chosen for truck and driver routes and how they can be synchronised.
The ILP formulations are defined in Section \ref{sec.ilp} and some families of valid inequalities in Section \ref{sec.valid.ineq}.
Section \ref{sec.computational.experiments} presents the computational experiments.
Finally, the conclusion and future work appear in Section \ref{sec.conclusions}.

\section{Problem description}\label{sec.problem.description}

This work considers a case study that involves the planning of a homogeneous fleet of trucks and drivers to transport pickup-and-delivery requests with multiple time windows over long distances.
The goal is to minimise the total travel cost of the trucks and taxis and the penalties for customer dissatisfaction caused by delivery delays.
More details are provided below and Fig. \ref{fig.notation.list} summarises the notation/definition list.

\begin{figure}
\begin{description}
%\deftitle{Scenario tree notation}
\item[$H \in \mathbb{N}_0$:] Number of days in the planning horizon.
\item[$I \in \mathbb{N}$:] Number of time instants per day, under certain discretization.
\item[{[$n$]}$ \doteq \{0,\ldots,n\}$:] Set of integers from 0 to $n$.
\item[$\mathcal{I} \doteq ${[$IH$]}:] Set of time instants in the planning horizon.
\item[$L,V,D,R$:] Set of locations, trucks, drivers, and requests, respectively.
\item[$l_v,l_d \in L$:] Start locations of $v \in V$ and $d \in D$, respectively.
\item[$L_V \doteq \{l_v : v \in V\}$:] Subset of truck start locations.
\item[$L_D \doteq \{l_d : d \in D\}$:] Subset of driver start locations.
\item[$l_r^{\textsc{p}},l_r^{\textsc{d}} \in L$:] Pickup and delivery locations of $r \in R$, respectively.
\item[$j_r^{\textsc{p}},j_r^{\textsc{d}} \in ${[$H-1$]}:] Pickup and delivery start days of $r \in R$, respectively.
\item[$a_r^{\textsc{p}},b_r^{\textsc{p}} \in ${[$I-1$]}:] Start and end times of the pickup time window of $r \in R$, respectively.
\item[$a_r^{\textsc{d}},b_r^{\textsc{d}} \in ${[$I-1$]}:] Start and end times of the delivery time window of $r \in R$, respectively.
\item[$s_r^{\textsc{p}},s_r^{\textsc{d}} \in \NN_0$:] Pickup and delivery service times of $r \in R$, respectively.
\item[$w_r \in \NN_0$:] Penalty for each day of delay in the delivery of $r \in R$.
\item[$\Day: \mathcal{I} \to ${[$H-1$]}:] Given $i \in \mathcal{I}$, $\Day(i)$ is the quotient of dividing $i$ by $I$.% i.e. the corresponding day in the planning horizon.}
\item[$\Time: \mathcal{I} \to ${[$I-1$]}:] Given $i \in \mathcal{I}$, $\Time(i)$ is the remainder of dividing $i$ by $I$.% i.e. the corresponding daily time instant.}
\item[$\VLEN,\TLEN: L \times L \to \NN$:] Truck and taxi travel time functions, respectively.
\item[$\VCOST,\TCOST: L \times L \to \NN$:] Truck and taxi travel cost functions, respectively.
\item[$\mathcal{I}_r^{\textsc{p}}, \mathcal{I}_r^{\textsc{d}} \subset \mathcal{I}$:] Subsets of instants where the pickup and delivery of $r \in R$ can begin, resp. %For example, if $r$ has a pickup time window as Fig. \ref{fig.example.time.window}(a), then $\mathcal{I}_r^{\textsc{p}} = \{8,\ldots,20\}\cup\{32,\ldots,44\}$.}
\end{description}
\caption{Notation/Definition list.}
\label{fig.notation.list}
\end{figure}

\begin{figure}
    \centering
    \includegraphics[width=\textwidth]{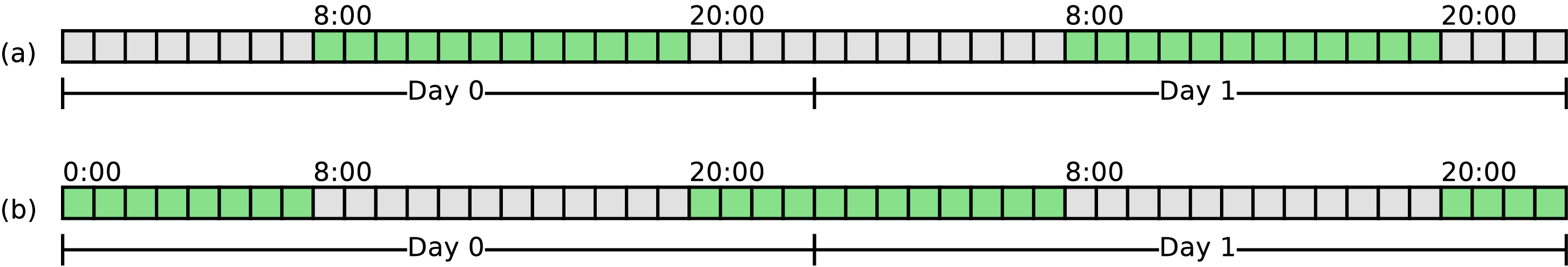}
    \caption{Examples of time windows in a two-day planning horizon.}
    \label{fig.example.time.window}
\end{figure}

Each request involves picking up (or loading) cargo in a location, transporting it to a different location, and delivering (or unloading) it there.
All requests must be delivered before the end of the planning horizon, typically 1 week.
The same truck that picks up a request must be in charge of delivering it, i.e. cargo transshipment is not allowed, and each truck can only transport a single request at a time.
Loading/unloading cargo into/from a truck takes 1 hour, known as service time, and the driver must remain in the truck in the meantime.
Each request has a pickup start day and a (hard) pickup time window indicating the day and time the loading is allowed to begin, but may end after hours.
Analogously, it also has a delivery start day and a (hard) delivery time window.
These time constraints represent the day from which the customer expects service and its hours of operation.
Each time window remains open from a start time to an end time (which can be on the next day as long as the 24-hour limit is not exceeded) and repeats this behaviour every day from its start day to the end of the planning horizon.
Fig. \ref{fig.example.time.window}a shows an example of a time window that opens from day 0 during the time interval [8:00, 20:00], where each box represents 1 hour.
Instead, the one in Fig. \ref{fig.example.time.window}b does so during [20:00, 8:00], but when this happens (closing on the next day), it must also open partially during [0:00, 8:00] on its start day (day 0) and during [20:00, 0:00] on the last day of the planning horizon (day 1).
Deliveries have no deadlines, but each request has a penalty for each day of delay from its delivery start day.
For example, for delivery time windows as Fig. \ref{fig.example.time.window}, there is no penalty on day 0, a 1-day penalty on day 1, and so on.

Each truck and driver has a start location, which represents their position at the beginning of planning.
They do not have to return there at the end of the planning horizon, but their last location becomes the start of the next planning.
There is a set of locations where trucks and drivers can stop, and they are not allowed to do so anywhere else, for example en route to some location.
This set contains all pickup and delivery locations of the requests and the start locations of the trucks and drivers.
The travel distance between locations is known and the travel time is constant throughout the day.

Drivers are the only ones authorised to drive the trucks, but they can also request an external taxi service to travel for an extra cost.
Therefore, crew scheduling becomes more flexible as drivers can take a taxi to a location other than their current position without using a truck.
There is an unlimited number of taxis available in every location and at every time, but each can carry one driver per ride.
A driver is considered to be working while travelling by vehicle (truck or taxi) and during service times and is supposed to be resting the rest of the time.
A rest begins when the driver gets out of a vehicle (truck or taxi) and ends as soon as he/she gets in the same or a different one.
Drivers must rest at least 12 hours in each 24-hour interval and must have at least 1 day off in each 7-day interval.
As mentioned, drivers can change trucks arbitrarily at any authorised location and at any time (before, during, or after service). 
For example, a driver who needs to rest can divert the truck to the nearest location, switch with another driver, and immediately continue transport without further interruptions. 
Drivers are not required to take a break after getting out of a vehicle, e.g. they can immediately switch to another.

Optionally, each truck can carry an extra driver as a passenger, who is also considered to be working.
Thus, crews can be composed of 1 or 2 drivers, and each one can be relieved independently.
Although this increases combinatorics even more, it may avoid taxis in certain situations.
For example, a driver can divert the truck to pick up a second driver and leave her/him in another location before continuing the route, instead of taking a taxi.
The choice between travelling as a passenger or taking a taxi will depend on the routes and the associated costs; e.g. it may be too expensive or impossible for a driver to travel as a passenger if all the trucks are far away at the time and need to travel long distances to her/his location.

Before proceeding to the mathematical models and more specific details, a toy instance of the problem is exemplified and two possible route plans are presented.
%to clarify the concepts defined so far and show the structure of the solutions.
%This example will be frequently referenced later to illustrate the construction of the digraphs.
%More challenging instances will be considered in the computational experiments.

\begin{example}\label{example}
Consider an instance with a planning horizon of one day discretised in 3-hour instants, i.e. $H = 1$ and $I = 8$.
There two locations $L=\{l_1,l_2\}$, two vehicles $V=\{v_1,v_2\}$, and two drivers $D=\{d_1,d_2\}$.
The start location of $v_1$ and $d_1$ is $l_1$ and that of $v_2$ and $d_2$ is $l_2$, i.e. $l_{v_1} = l_{d_1} = l_1$ and $l_{v_2} = l_{d_2} = l_2$.
The travel time between both locations is one instant, i.e. $\VLEN(l_1,l_2) = \VLEN(l_2,l_1) = \TLEN(l_1,l_2) = \TLEN(l_2,l_1) = 1$.
There are two requests $R = \{r_1,r_2\}$, whose attributes appear in Table \ref{tab.request.attributes}.

\begin{table}
\caption{Request attributes in Example \ref{example}}%
\label{tab.request.attributes}
\begin{tabular}{@{}@{\extracolsep{\fill}}llllllllll@{}}
\hline
Request & \multicolumn{4}{l}{Pickup} &  \multicolumn{4}{l}{Delivery} \\ [-6pt]
& \multicolumn{4}{l}{\hrulefill} & \multicolumn{4}{l@{}}{\hrulefill} \\
~ & Location & Start day & Time w. & Service t. & Location & Start day & Time w. & Service t. \\
$r \in R$ & $l^\textsc{p}_r$ & $j^\textsc{p}_r$ & $[a^\textsc{p}_r, b^\textsc{p}_r]$ & $s^\textsc{p}_r$ & $l^\textsc{d}_r$ & $j^\textsc{d}_r$ & $[a^\textsc{d}_r, b^\textsc{d}_r]$ &  $s^\textsc{d}_r$ \\
\hline
$r_1$ & $l_1$ & 0 & $[0,2]$ & 1 & $l_2$ & 0 & $[6,2]$ & 1\\
$r_2$ & $l_2$ & 0 & $[3,5]$ & 1 & $l_1$ & 0 & $[3,5]$ & 1\\
\hline
\end{tabular}
\end{table}

Fig. \ref{fig.example}a depicts possible routes for the trucks and drivers.
Recall that the delivery of $r_1$ is allowed to begin at 6:00 since its time window opens during [0,2] = [0:00,6:00] and [6,8] = [18:00,0:00] on day 0, as explained in Fig. \ref{fig.example.time.window}b.
Fig. \ref{fig.example}b depicts a better route for $d_2$, where he/she travels as a passenger. %in a truck instead of taking a taxi.

\begin{figure}
\centering
\includegraphics[width=\textwidth]{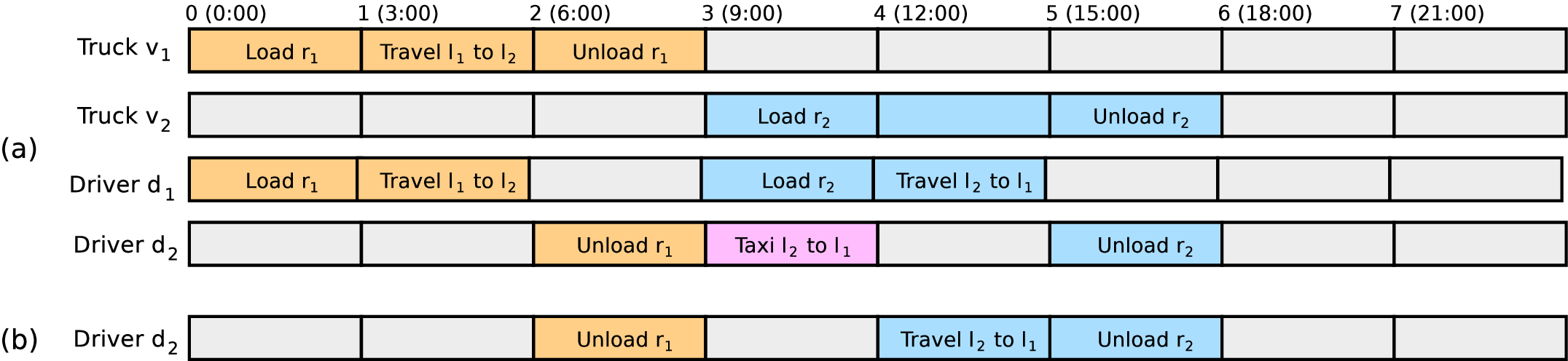}
\caption{Routes of the Example \ref{example}. The tasks performed by truck $v_1$ are orange and those by $v_2$ are blue.}
\label{fig.example}  
\end{figure}

\end{example}

\section{Route representation}\label{sec.routes.representation}

This section defines digraphs that represent truck and driver routes as directed paths.
The arcs specify the actions performed over time, which can be resting, travelling, loading, or unloading.
However, there are some categories of directed paths that do not correspond to valid routes and should be prohibited.
In the end, the conditions that routes must meet to be synchronised in space and time are identified.

\subsection{Truck routes}\label{sec.models.truck.routes}

\vspace{-5pt}
\paragraph{Location--Time digraph}

Let $G_{\textsc{lt}} \doteq (N,E,c)$ be a weighted multidigraph constructed as follows.
For each $l \in L$ and $i \in \mathcal{I}$, there is a node $n_{l,i}$ in $N$, representing that a truck is in location $l$ at instant $i$.
There are also two distinguished nodes in $N$, the source $n_s$ and the sink $\tilde n_{s}$.
The multiset $E$ is the union of the following sets.

\begin{figure}
\centering

\begin{tikzpicture}

% Source and sink nodes
\node[fill,circle,inner sep=2pt, label=west:$n_s$] (0) at (0,0.5) {};
\node[fill,circle,inner sep=2pt, label=east:$\tilde n_s$] (1) at (10,0.5) {};

% Rest of the nodes
\foreach \i in {1,...,9}
{
 	\foreach \l [evaluate={\y = \l/2}] in {0,2}
     {
 		\node[fill,circle,inner sep=2pt] (\l\i) at (\i,\y) {};
 	}       
}

% Node labels
\node at (0.5,-0.2) {$n_{l_2,0}$};
\node at (0.5,1.2) {$n_{l_1,0}$}; 
\node at (9.65,-0.2) {$n_{l_2,|\mathcal{I}|}$};
\node at (9.65,1.1) {$n_{l_1,|\mathcal{I}|}$}; 

% Truck arcs
\foreach \i [evaluate={\j = int(\i+1)}] in {1,2,3,4,6,7,8}
{
    \draw[->,line width=0.5pt] (0\i) to (2\j);
}
\foreach \i [evaluate={\j = int(\i+1)}] in {1,3,4,...,8}
{
    \draw[->,line width=0.5pt] (2\i) to (0\j);
}

% Rest arcs
\foreach \i [evaluate={\j = int(\i+1)}] in {1,...,6}
{
    \draw[->,line width=0.5pt] (2\i) to (2\j);
}

% Pickup and delivery arcs of request 1
\foreach \i [evaluate={\j = int(\i+1)}, evaluate={\x = \i + 0.5}] in {2,3}
{
    \draw[->,line width=0.5pt, bend left = 90, looseness = 1.5] (2\i) to (2\j) node [fill=white, inner sep=1pt] at (\x,1.5) {\scriptsize ${\uparrow}r_1$};
}
\foreach \i [evaluate={\j = int(\i+1)}, evaluate={\x = \i + 0.5}] in {1,2,7,8}
{
    \draw[->,line width=0.5pt, bend right = 90, looseness = 1.5] (0\i) to (0\j) node [fill=white, inner sep=1pt] at (\x,-0.55) {\scriptsize ${\downarrow}r_1$};
}

% Pickup and delivery arcs of request 2
\foreach \i [evaluate={\j = int(\i+1)}, evaluate={\x = \i + 0.5}] in {5,6}
{
    \draw[->,line width=0.5pt, bend right = 90, looseness = 1.5] (0\i) to (0\j) node [fill=white, inner sep=1pt] at (\x,-0.55) {\scriptsize ${\uparrow}r_2$};
}
\foreach \i [evaluate={\j = int(\i+1)}, evaluate={\x = \i + 0.5}] in {4,5}
{
    \draw[->,line width=0.5pt, bend left = 90, looseness = 1.5] (2\i) to (2\j) node [fill=white, inner sep=1pt] at (\x,1.5) {\scriptsize ${\downarrow}r_2$};
}

\draw[line width=0.5mm, orange, ->] (0) to (21);
\draw[->,line width=0.5mm, bend left = 90, looseness = 1.5, orange] (21) to (22) node [fill=white, inner sep=1pt] at (1.5,1.5) {\scriptsize ${\uparrow}r_1$};
\draw[line width=0.5mm, orange, ->] (22) to (03);
\draw[->,line width=0.5mm, bend right = 90, looseness = 1.5, orange] (03) to (04) node [fill=white, inner sep=1pt] at (3.5,-0.55) {\scriptsize ${\downarrow}r_1$};
\draw[line width=0.5mm, orange, ->] (04) to (05);
\draw[line width=0.5mm, orange, ->] (05) to (06);
\draw[line width=0.5mm, orange, ->] (06) to (07);
\draw[line width=0.5mm, orange, ->] (07) to (08);
\draw[line width=0.5mm, orange, ->] (08) to (09);
\draw[line width=0.5mm, orange, ->] (09) to (1);

\draw[line width=0.5mm, cyan, ->] (0) to (01);
\draw[line width=0.5mm, cyan, ->] (01) to (02);
\draw[line width=0.5mm, cyan, ->] (02) to (03);
\draw[line width=0.5mm, cyan, ->] (03) to (04);
\draw[->,line width=0.5mm, bend right = 90, looseness = 1.5, cyan] (04) to (05) node [fill=white, inner sep=1pt] at (4.5,-0.55) {\scriptsize ${\uparrow}r_2$};
\draw[line width=0.5mm, cyan, ->] (05) to (26);
\draw[->,line width=0.5mm, bend left = 90, looseness = 1.5, cyan] (26) to (27) node [fill=white, inner sep=1pt] at (6.5,1.5) {\scriptsize ${\downarrow}r_2$};
\draw[line width=0.5mm, cyan, ->] (27) to (28);
\draw[line width=0.5mm, cyan, ->] (28) to (29);
\draw[line width=0.5mm, cyan, ->] (29) to (1);

\end{tikzpicture}

\caption{Digraph $G_{\textsc{lt}}$ for Example \ref{example}.}
\label{fig.digraph.LT}
\end{figure}
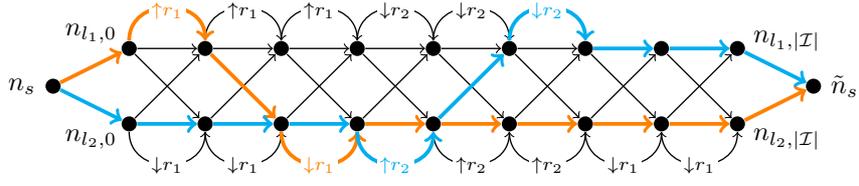

\begin{itemize}

    \item $E^{\textsc{rest}}$ has an arc $(n_{l,i}, n_{l,i+1})$ for each $l \in L$ and $i \in [|\mathcal{I}|-1]$, representing a rest of one time instant.

    \item $E^{\textsc{trip}}$ has an arc $(n_{l_1,i},n_{l_2,i + \Delta})$ for each pair $(l_1,l_2) \in L \times L$ of adjacent locations and $i \in [|\mathcal{I}| - \Delta]$ with  $\Delta = \VLEN(l_1,l_2)$, representing that a truck travels.
    
    \item For each $r \in R$, $E^{\textsc{p},r}$ has an arc $(n_{l^\textsc{p}_r,i},n_{l^\textsc{p}_r,i+s^\textsc{p}_r})$ for each $i \in \mathcal{I}^\textsc{p}_r$, representing that a truck loads $r$.

    \item For each $r \in R$, $E^{\textsc{d},r}$ has an arc $(n_{l^\textsc{d}_r,i},n_{l^\textsc{d}_r,i+s^\textsc{d}_r})$ for each $i \in \mathcal{I}^\textsc{d}_r$, representing that a truck unloads $r$.
    
    \item $E^{\textsc{sour}}$ has an arc $(n_s,n_{l,0})$ for each $l \in L_V$ and $E^{\textsc{sink}}$ has an arc $(n_{l,|\mathcal{I}|}, \tilde n_s)$ for each $l \in L$.

\end{itemize}
The weight of an arc $e = (n_{l_1,i},n_{l_2,i + \Delta}) \in E^{\textsc{trip}}$ is the truck travel cost between $l_1$ and $l_2$, i.e. $c(e) = \VCOST(l_1,l_2)$, the weight of $e' =  (n_{l^\textsc{d}_r,i},n_{l^\textsc{d}_r,i+s^\textsc{d}_r}) \in E^{\textsc{d},r}$ is the product between the penalty of $r$ and the number of days of delay, i.e. $c(e') = w_r(\Day(i) - j^{\textsc{d}}_r)$, and the remaining are zero-weighted arcs.
Fig. \ref{fig.digraph.LT} shows the resulting construction for the instance of Example \ref{example}.
The position of the nodes on the vertical and horizontal axis determines the location (first subscript) and the time instant (second subscript), respectively.
The arcs in $E^{\textsc{rest}}$ are parallel to the horizontal axis and those in $E^{\textsc{trip}}$ are oblique and not incident on the source/sink.
For each $r \in R$, the ones in $E^{\textsc{p},r}$ and $E^{\textsc{d},r}$ are curved with the label ``${\uparrow}r$'' and ``${\downarrow}r$'', respectively, to distinguish them from other parallel arcs.
The orange $(n_s,\tilde n_s)$-directed path corresponds to the route of the truck $v_1$ and the blue one to $v_2$.

When there are ambiguities, the node and arc sets of $G_{\textsc{lt}}$ will also have the subscript ``$\textsc{lt}$'', e.g. $N_{\textsc{lt}}$, $E_{\textsc{lt}}$, $E_{\textsc{lt}}^{\textsc{rest}}$, etc., as well as the weight vector $c_{\textsc{lt}}$.
To avoid referring to particular requests, the multisets $E^{\textsc{p}}$ and $E^{\textsc{d}}$ will denote ${\textstyle \bigcup_{r \in R}E^{\textsc{p},r}}$ and ${\textstyle \bigcup_{r \in R}E^{\textsc{d},r}}$, respectively.
By construction, $G_{\textsc{lt}}$ is an acyclic digraph where every truck route corresponds to a $(n_s,\tilde n_s)$-directed path. 
However, the converse is not necessarily true.
Let $P$ be a $(n_s,\tilde n_s)$-directed path, four categories of forbidden paths can be identified based on the requirements of the case study.

\begin{figure}
\centering

\begin{tikzpicture}

\node at (-1.5,0.5) {(a)};

% Source and sink nodes
\node[fill,circle,inner sep=2pt, label=west:$n_s$] (0) at (0,0.5) {};
\node[fill,circle,inner sep=2pt, label=east:$\tilde n_s$] (1) at (10,0.5) {};

% Rest of the nodes
\foreach \i in {1,...,9}
{
    \node[fill,circle,inner sep=2pt] (0\i) at (\i,0) {};
    \node[fill,circle,inner sep=2pt] (2\i) at (\i,1) {};     
}

% Node labels
\node at (0.5,-0.2) {$n_{l_2,0}$};
\node at (0.5,1.2) {$n_{l_1,0}$}; 
\node at (9.65,-0.2) {$n_{l_2,|\mathcal{I}|}$};
\node at (9.65,1.1) {$n_{l_1,|\mathcal{I}|}$}; 

% Truck arcs
\foreach \i [evaluate={\j = int(\i+1)}] in {1,...,8}
{
    \draw[->,line width=0.5pt] (0\i) to (2\j);
    \draw[->,line width=0.5pt] (2\i) to (0\j);
}

% Rest arcs
\foreach \i [evaluate={\j = int(\i+1)}] in {4,7}
{
    \draw[->,line width=0.5pt] (0\i) to (0\j);
}
\foreach \i [evaluate={\j = int(\i+1)}] in {1,4}
{
    \draw[->,line width=0.5pt] (2\i) to (2\j);
}

% Pickup and delivery arcs of request 1
\foreach \i [evaluate={\j = int(\i+1)}, evaluate={\x = \i + 0.5}] in {2,3}
{
    \draw[->,line width=0.5pt, bend left = 90, looseness = 1.5] (2\i) to (2\j) node [fill=white, inner sep=1pt] at (\x,1.5) {\scriptsize ${\uparrow}r_1$};
}
\foreach \i [evaluate={\j = int(\i+1)}, evaluate={\x = \i + 0.5}] in {1,2,3,8}
{
    \draw[->,line width=0.5pt, bend right = 90, looseness = 1.5] (0\i) to (0\j) node [fill=white, inner sep=1pt] at (\x,-0.55) {\scriptsize ${\downarrow}r_1$};
}

% Pickup and delivery arcs of request 2
\foreach \i [evaluate={\j = int(\i+1)}, evaluate={\x = \i + 0.5}] in {5,6}
{
    \draw[->,line width=0.5pt, bend right = 90, looseness = 1.5] (0\i) to (0\j) node [fill=white, inner sep=1pt] at (\x,-0.55) {\scriptsize ${\uparrow}r_2$};
}
\foreach \i [evaluate={\j = int(\i+1)}, evaluate={\x = \i + 0.5}] in {5,6}
{
    \draw[->,line width=0.5pt, bend left = 90, looseness = 1.5] (2\i) to (2\j) node [fill=white, inner sep=1pt] at (\x,1.5) {\scriptsize ${\downarrow}r_2$};
}

\draw[line width=0.5mm, orange, ->] (0) to (21);
\draw[->,line width=0.5mm, bend left = 90, looseness = 1.5, orange] (21) to (22) node [fill=white, inner sep=1pt] at (1.5,1.5) {\scriptsize ${\uparrow}r_1$};
\draw[line width=0.5mm, orange, ->] (22) to (23);
\draw[line width=0.5mm, orange, ->] (23) to (24);
\draw[->,line width=0.5mm, bend left = 90, looseness = 1.5, orange] (24) to (25) node [fill=white, inner sep=1pt] at (4.5,1.5) {\scriptsize ${\downarrow}r_2$};
\draw[line width=0.5mm, orange, ->] (25) to (26);
\draw[line width=0.5mm, orange, ->] (26) to (27);
\draw[line width=0.5mm, orange, ->] (27) to (28);
\draw[line width=0.5mm, orange, ->] (28) to (29);
\draw[line width=0.5mm, orange, ->] (29) to (1);

\draw[line width=0.5mm, cyan, ->] (0) to (01);
\draw[line width=0.5mm, cyan, ->] (01) to (02);
\draw[line width=0.5mm, cyan, ->] (02) to (03);
\draw[line width=0.5mm, cyan, ->] (03) to (04);
\draw[->,line width=0.5mm, bend right = 90, looseness = 1.5, cyan] (04) to (05) node [fill=white, inner sep=1pt] at (4.5,-0.55) {\scriptsize ${\uparrow}r_2$};
\draw[line width=0.5mm, cyan, ->] (05) to (06);
\draw[line width=0.5mm, cyan, ->] (06) to (07);
\draw[->,line width=0.5mm, bend right = 90, looseness = 1.5, cyan] (07) to (08) node [fill=white, inner sep=1pt] at (7.5,-0.55) {\scriptsize ${\downarrow}r_1$};
\draw[line width=0.5mm, cyan, ->] (08) to (09);
\draw[line width=0.5mm, cyan, ->] (09) to (1);

\end{tikzpicture}

\medskip

\begin{tikzpicture}

\node at (-1.5,0.5) {(b)};

% Source and sink nodes
\node[fill,circle,inner sep=2pt, label=west:$n_s$] (0) at (0,0.5) {};
\node[fill,circle,inner sep=2pt, label=east:$\tilde n_s$] (1) at (10,0.5) {};

% Rest of the nodes
\foreach \i in {1,...,9}
{
    \node[fill,circle,inner sep=2pt] (0\i) at (\i,0) {};
    \node[fill,circle,inner sep=2pt] (2\i) at (\i,1) {};     
}

% Rest arcs
\foreach \i [evaluate={\j = int(\i+1)}] in {1,...,6,8}
{
    \draw[->,line width=0.5pt] (0\i) to (0\j);
}
\foreach \i [evaluate={\j = int(\i+1)}] in {2,4,5,...,8}
{
    \draw[->,line width=0.5pt] (2\i) to (2\j);
}

% Node labels
\node at (0.5,-0.2) {$n_{l_2,0}$};
\node at (0.5,1.2) {$n_{l_1,0}$}; 
\node at (9.65,-0.2) {$n_{l_2,|\mathcal{I}|}$};
\node at (9.65,1.1) {$n_{l_1,|\mathcal{I}|}$};  

% Source and sink arcs
\foreach \l in {0}
{
    \draw[->,line width=0.5pt] (0) to (\l1);   
} 
\foreach \l in {2}
{
    \draw[->,line width=0.5pt] (\l9) to (1);    
}

% Truck arcs
\foreach \i [evaluate={\j = int(\i+1)}] in {1,...,8}
{
    \foreach \l [evaluate={\m = int(\l+2)}] in {0}
    {
    	\draw[->,line width=0.5pt] (\l\i) to (\m\j);
    	\draw[->,line width=0.5pt] (\m\i) to (\l\j);  %      
    }
}

% Pickup and delivery arcs of request 1
\foreach \i [evaluate={\j = int(\i+1)}, evaluate={\x = \i + 0.5}] in {1,3}
{
    \draw[->,line width=0.5pt, bend left = 90, looseness = 1.5] (2\i) to (2\j) node [fill=white, inner sep=1pt] at (\x,1.5) {\scriptsize ${\uparrow}r_1$};
}
\foreach \i [evaluate={\j = int(\i+1)}, evaluate={\x = \i + 0.5}] in {1,2,3,7}
{
    \draw[->,line width=0.5pt, bend right = 90, looseness = 1.5] (0\i) to (0\j) node [fill=white, inner sep=1pt] at (\x,-0.55) {\scriptsize ${\downarrow}r_1$};
}

% Pickup and delivery arcs of request 2
\foreach \i [evaluate={\j = int(\i+1)}, evaluate={\x = \i + 0.5}] in {4,5}
{
    \draw[->,line width=0.5pt, bend right = 90, looseness = 1.5] (0\i) to (0\j) node [fill=white, inner sep=1pt] at (\x,-0.55) {\scriptsize ${\uparrow}r_2$};
}
\foreach \i [evaluate={\j = int(\i+1)}, evaluate={\x = \i + 0.5}] in {5,6}
{
    \draw[->,line width=0.5pt, bend left = 90, looseness = 1.5] (2\i) to (2\j) node [fill=white, inner sep=1pt] at (\x,1.5) {\scriptsize ${\downarrow}r_2$};
}

\draw[line width=0.5mm, orange, ->] (0) to (21);
\draw[line width=0.5mm, orange, ->] (21) to (22);
\draw[->,line width=0.5mm, bend left = 90, looseness = 1.5, orange] (22) to (23) node [fill=white, inner sep=1pt] at (2.5,1.5) {\scriptsize ${\uparrow}r_1$};
\draw[line width=0.5mm, orange, ->] (23) to (24);
\draw[->,line width=0.5mm, bend left = 90, looseness = 1.5, orange] (24) to (25) node [fill=white, inner sep=1pt] at (4.5,1.5) {\scriptsize ${\downarrow}r_2$};
\draw[line width=0.5mm, orange, ->] (25) to (06);
\draw[->,line width=0.5mm, bend right = 90, looseness = 1.5, orange] (06) to (07) node [fill=white, inner sep=1pt] at (6.5,-0.55) {\scriptsize ${\uparrow}r_2$};
\draw[line width=0.5mm, orange, ->] (07) to (08);
\draw[->,line width=0.5mm, bend right = 90, looseness = 1.5, orange] (08) to (09) node [fill=white, inner sep=1pt] at (8.5,-0.55) {\scriptsize ${\downarrow}r_1$};
\draw[line width=0.5mm, orange, ->] (09) to (1);

\end{tikzpicture}

\medskip

\begin{tikzpicture}

\node at (-1.5,0.5) {(c)};

% Source and sink nodes
\node[fill,circle,inner sep=2pt, label=west:$n_s$] (0) at (0,0.5) {};
\node[fill,circle,inner sep=2pt, label=east:$\tilde n_s$] (1) at (10,0.5) {};

% Rest of the nodes
\foreach \i in {1,...,9}
{
    \node[fill,circle,inner sep=2pt] (0\i) at (\i,0) {};
    \node[fill,circle,inner sep=2pt] (2\i) at (\i,1) {};     
}

% Node labels
\node at (0.5,-0.2) {$n_{l_2,0}$};
\node at (0.5,1.2) {$n_{l_1,0}$}; 
\node at (9.65,-0.2) {$n_{l_2,|\mathcal{I}|}$};
\node at (9.65,1.1) {$n_{l_1,|\mathcal{I}|}$};  

% Source and sink arcs
\foreach \l in {0}
{
    \draw[->,line width=0.5pt] (0) to (\l1);   
} 
\foreach \l in {2}
{
    \draw[->,line width=0.5pt] (\l9) to (1);    
} 

% Truck arcs
\foreach \i [evaluate={\j = int(\i+1)}] in {1,2,3,4,6,7,8}
{
    \draw[->,line width=0.5pt] (0\i) to (2\j);
}
\foreach \i [evaluate={\j = int(\i+1)}] in {1,2,4,5,6,8}
{
    \draw[->,line width=0.5pt] (2\i) to (0\j);
}

% Rest arcs
\foreach \i [evaluate={\j = int(\i+1)}] in {1,...,8}
{
	\foreach \l in {0,2}
    {
    	\draw[->,line width=0.5pt] (\l\i) to (\l\j);
    }
}

% Pickup and delivery arcs of request 1
\foreach \i [evaluate={\j = int(\i+1)}, evaluate={\x = \i + 0.5}] in {2,3}
{
    \draw[->,line width=0.5pt, bend left = 90, looseness = 1.5] (2\i) to (2\j) node [fill=white, inner sep=1pt] at (\x,1.5) {\scriptsize ${\uparrow}r_1$};
}
\foreach \i [evaluate={\j = int(\i+1)}, evaluate={\x = \i + 0.5}] in {1,2,3,7}
{
    \draw[->,line width=0.5pt, bend right = 90, looseness = 1.5] (0\i) to (0\j) node [fill=white, inner sep=1pt] at (\x,-0.55) {\scriptsize ${\downarrow}r_1$};
}

% Pickup and delivery arcs of request 2
\foreach \i [evaluate={\j = int(\i+1)}, evaluate={\x = \i + 0.5}] in {5,6}
{
    \draw[->,line width=0.5pt, bend right = 90, looseness = 1.5] (0\i) to (0\j) node [fill=white, inner sep=1pt] at (\x,-0.55) {\scriptsize ${\uparrow}r_2$};
}
\foreach \i [evaluate={\j = int(\i+1)}, evaluate={\x = \i + 0.5}] in {4,5}
{
    \draw[->,line width=0.5pt, bend left = 90, looseness = 1.5] (2\i) to (2\j) node [fill=white, inner sep=1pt] at (\x,1.5) {\scriptsize ${\downarrow}r_2$};
}

\draw[line width=0.5mm, orange, ->] (0) to (21);
\draw[->,line width=0.5mm, bend left = 90, looseness = 1.5, orange] (21) to (22) node [fill=white, inner sep=1pt] at (1.5,1.5) {\scriptsize ${\uparrow}r_1$};
\draw[line width=0.5mm, orange, ->] (22) to (23);
\draw[line width=0.5mm, orange, ->] (23) to (04);
\draw[->,line width=0.5mm, bend right = 90, looseness = 1.5, orange] (04) to (05) node [fill=white, inner sep=1pt] at (4.5,-0.55) {\scriptsize ${\uparrow}r_2$};
\draw[line width=0.5mm, orange, ->] (05) to (26);
\draw[->,line width=0.5mm, bend left = 90, looseness = 1.5, orange] (26) to (27) node [fill=white, inner sep=1pt] at (6.5,1.5) {\scriptsize ${\downarrow}r_2$};
\draw[line width=0.5mm, orange, ->] (27) to (08);
\draw[->,line width=0.5mm, bend right = 90, looseness = 1.5, orange] (08) to (09) node [fill=white, inner sep=1pt] at (8.5,-0.55) {\scriptsize ${\downarrow}r_1$}; 
\draw[line width=0.5mm, orange, ->] (09) to (1);

\end{tikzpicture}

\caption{Examples of forbidden paths in $G_{\textsc{lt}}$.}
\label{fig.lt.forbidden.paths}
\end{figure}
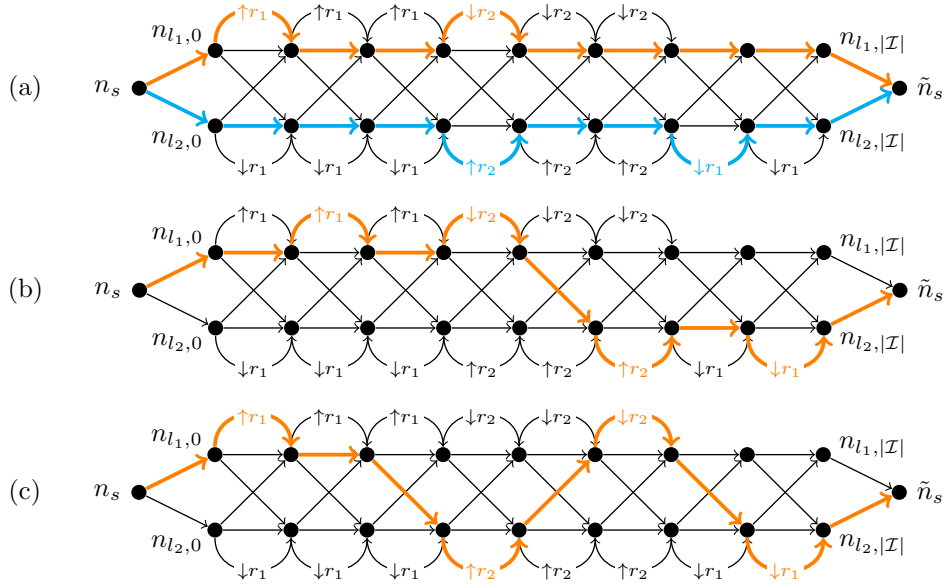

\begin{itemize}
\item \textit{Repeated services}. 
Since each request must be picked up and delivered exactly once, then for all $r \in R$, $P$ cannot have more than one arc of $E^{\textsc{p},r}$ and $E^{\textsc{d},r}$. 
%the number of arcs with labels ${\uparrow}r$ and ${\downarrow}r$ in $P$ cannot be greater than 1, i.e. 
%$|P \cap E^{\textsc{p},r}| \leq 1$ and $|P \cap E^{\textsc{d},r}| \leq 1$.
\item \textit{Unpaired services}. 
Since each request must be fully served by the same truck, then for all $r \in R$, the total number of arcs of $E^{\textsc{p},r}$ and $E^{\textsc{d},r}$ in $P$ cannot differ.
%the number of arcs with labels ${\uparrow}r$ and ${\downarrow}r$ in $P$ cannot differ, i.e. 
%$|P \cap E^{\textsc{p},r}| = |P \cap E^{\textsc{d},r}|$.
Two of these forbidden paths for the instance in Example \ref{example} are depicted in Fig. \ref{fig.lt.forbidden.paths}(a).
\item \textit{Disordered services}. 
Since each request must be picked up before being delivered, 
then for all $r \in R$ and $P'$ subdirected path of $P$ from $n_s$ such that its last arc belongs to $E^{\textsc{d},r}$, the total number of arcs of $E^{\textsc{p},r}$ in $P'$ cannot be zero.
% then for all $r \in R$ and $e \in P \cap E^{\textsc{d},r}$, the total number of arcs of $E^{\textsc{p},r}$ in the subdirected path of $P$ from $n_s$ to the tail of $e$ cannot be zero.
% then for all $P'$ subdirected path of $P$ from $n_s$, the difference between the total number of arcs of $E^{\textsc{d},r}$ and $E^{\textsc{p},r}$ in $P'$ cannot be greater than 0.
%then for all $r \in R$, and $e \in P \cap E^{\textsc{d},r}$, then the number of arcs in $P \cap E^{\textsc{p},r}$ before $e$ cannot be zero.
Fig. \ref{fig.lt.forbidden.paths}(b) shows an example of such a forbidden path that also has no unpaired services.
\item \textit{Excess capacity.} 
Since each truck can service only one request at a time, then for all $P'$ subdirected path of $P$ from $n_s$ such that its last arc belongs to $E^{\textsc{p}}$, the difference between the total number of arcs of $E^{\textsc{p}}$ and $E^{\textsc{d}}$ in $P'$ cannot be greater than one.
%for all $r_1,r_2 \in R$, $e_1 \in P \cap E^{\textsc{p},r_1}$ and $e_2 \in P \cap E^{\textsc{p},r_2}$, the number of arcs in $P \cap (E^{\textsc{d},r_1} \cup E^{\textsc{d},r_2})$ between $e_1$ and $e_2$ cannot be zero.
Fig. \ref{fig.lt.forbidden.paths}(c) shows an example of such a forbidden path that also has no unpaired or disordered services.
\end{itemize}

In the remainder of the subsection, digraphs with greater structure are proposed that prevent the existence of certain types of forbidden paths.

\vspace{-5pt}
\paragraph{Location--Time--Cargo digraph}

Let $G_{\textsc{ltc}} \doteq (N,E,c)$ be a weighted multidigraph constructed as follows.
For each $l \in L$, $i \in \mathcal{I}$, and $q \in \{\vacio,\lleno\}$, there is a node $n_{l,i}^q$ in $N$, representing that the truck is empty in location $l$ at instant $i$ if $q=\vacio$ and is carrying cargo if $q=\lleno$.
There are also two distinguished nodes in $N$, the source $n_s$ and the sink $\tilde n_{s}$.
The multiset $E$ follows the construction of $E_{\textsc{lt}}$, taking into account that the pickup and delivery arcs are the only ones that modify cargo and that trucks must begin and end their route empty.
Formally, $E$ is the union of the following sets.

\begin{figure}
\centering
\begin{tikzpicture}

% Source and sink nodes
\node[fill,circle,inner sep=2pt, label=west:$n_s$] (0) at (0,0.5) {};
\node[fill,circle,inner sep=2pt, label=east:$\tilde n_s$] (1) at (10,0.5) {};

% Rest of the nodes
\foreach \i in {1,...,9}
{
    \node[fill,circle,inner sep=2pt] (2\i) at (\i,0) {};
}
\foreach \i in {1,...,9}
{
    \node[fill,circle,inner sep=2pt] (3\i) at (\i,1) {};
}
\foreach \i in {1,2,9}
{
    \node[draw,circle,inner sep=2pt,black!30,densely dashed] (4\i) at (\i,2) {};
}
\foreach \i in {3,4,...,8}
{
    \node[fill,circle,inner sep=2pt] (4\i) at (\i,2) {};
}
\foreach \i in {1,8,9}
{
    \node[draw,circle,inner sep=2pt,black!30,densely dashed] (5\i) at (\i,3) {};
}
\foreach \i in {2,3,...,7}
{
    \node[fill,circle,inner sep=2pt] (5\i) at (\i,3) {};
}

% Node labels
\node at (0.5,-0.1) {$n^{\vacio}_{l_2,0}$};
\node at (0.5,1.2) {$n^{\vacio}_{l_1,0}$}; 
\node at (0.5,1.9) {$n^{\lleno}_{l_2,0}$};
\node at (0.5,3) {$n^{\lleno}_{l_1,0}$}; 
\node at (9.7,-0.1) {$n^{\vacio}_{l_2,|\mathcal{I}|}$};
\node at (9.7,1.2) {$n^{\vacio}_{l_1,|\mathcal{I}|}$}; 
\node at (9.7,1.9) {$n^{\lleno}_{l_2,|\mathcal{I}|}$};
\node at (9.7,3) {$n^{\lleno}_{l_1,|\mathcal{I}|}$}; 

% Truck arcs
\foreach \i [evaluate={\j = int(\i+1)}] in {1,...,8}
{
    \draw[->,line width=0.5pt] (2\i) to (3\j);
    \draw[->,line width=0.5pt] (3\i) to (2\j);
}
\foreach \i [evaluate={\j = int(\i+1)}] in {1,2,7,8}
{
    \draw[->,line width=0.5pt, black!30, dashed] (4\i) to (5\j);
}
\foreach \i [evaluate={\j = int(\i+1)}] in {3,4,6}
{
    \draw[->,line width=0.5pt] (4\i) to (5\j);
}
\foreach \i [evaluate={\j = int(\i+1)}] in {1,8}
{
    \draw[->,line width=0.5pt, black!30, dashed] (5\i) to (4\j);
}
\foreach \i [evaluate={\j = int(\i+1)}] in {3,4,5,6,7}
{
    \draw[->,line width=0.5pt] (5\i) to (4\j);
}

% Rest arcs
\foreach \i [evaluate={\j = int(\i+1)}] in {1,...,6}
{
    \draw[->,line width=0.5pt] (3\i) to (3\j);
}
\foreach \i [evaluate={\j = int(\i+1)}] in {1,2,8}
{
    \draw[->,line width=0.5pt,black!30,dashed] (4\i) to (4\j);
}
\foreach \i [evaluate={\j = int(\i+1)}] in {3,4,...,7}
{
    \draw[->,line width=0.5pt] (4\i) to (4\j);
}
\foreach \i [evaluate={\j = int(\i+1)}] in {1,7,8}
{
    \draw[->,line width=0.5pt,black!30,dashed] (5\i) to (5\j);
}
\foreach \i [evaluate={\j = int(\i+1)}] in {2,3,...,6}
{
    \draw[->,line width=0.5pt] (5\i) to (5\j);
}

% Pickup and delivery arcs of request 1
\foreach \i [evaluate={\j = int(\i+1)}, evaluate={\x = \i + 0.4}] in {2,3}
{
    \draw[->,line width=0.5pt] (3\i) to (5\j) node[fill=white,inner sep=0.5pt] at (\x,1.8) {\scriptsize ${\uparrow}r_1$};
}
\foreach \i [evaluate={\j = int(\i+1)}, evaluate={\x = \i + 0.4}] in {7,8}
{
    \draw[->,line width=0.5pt] (4\i) to (2\j) node[fill=white,inner sep=0.5pt] at (\x,1.2) {\scriptsize ${\downarrow}r_1$};
}
\foreach \i [evaluate={\j = int(\i+1)}, evaluate={\x = \i + 0.4}] in {1,2}
{
    \draw[->,line width=0.5pt,black!30,dashed] (4\i) to (2\j) node[fill=white,inner sep=0.5pt] at (\x,1.2) {\scriptsize ${\downarrow}r_1$};
}

% Pickup and delivery arcs of request 2
\foreach \i [evaluate={\j = int(\i+1)}, evaluate={\x = \i + 0.6}] in {5,6}
{
    \draw[->,line width=0.5pt] (2\i) to (4\j) node [fill=white,inner sep=0.5pt] at (\x,1.2) {\scriptsize ${\uparrow}r_2$};
}
\foreach \i [evaluate={\j = int(\i+1)}, evaluate={\x = \i + 0.6}] in {4,5}
{
    \draw[->,line width=0.5pt] (5\i) to (3\j) node[fill=white,inner sep=0.5pt] at (\x,1.8) {\scriptsize ${\downarrow}r_2$};
}

\draw[line width=0.5mm, orange, ->] (0) to (31);
\draw[->,line width=0.5mm, orange] (31) to (52) node [fill=white,inner sep=0.5pt] at (1.4,1.8) {\scriptsize ${\uparrow}r_1$};
\draw[line width=0.5mm, orange, ->] (52) to (43);
\draw[->,line width=0.5mm, orange] (43) to (24) node[fill=white,inner sep=0.5pt] at (3.4,1.2) {\scriptsize ${\downarrow}r_1$};
\draw[line width=0.5mm, orange, ->] (24) to (25);
\draw[line width=0.5mm, orange, ->] (25) to (26);
\draw[line width=0.5mm, orange, ->] (26) to (27);
\draw[line width=0.5mm, orange, ->] (27) to (28);
\draw[line width=0.5mm, orange, ->] (28) to (29);
\draw[line width=0.5mm, orange, ->] (29) to (1);

\draw[line width=0.5mm, cyan, ->] (0) to (21);
\draw[line width=0.5mm, cyan, ->] (21) to (22);
\draw[line width=0.5mm, cyan, ->] (22) to (23);
\draw[line width=0.5mm, cyan, ->] (23) to (24);
\draw[->,line width=0.5mm, cyan] (24) to (45) node [fill=white,inner sep=0.5pt] at (4.6,1.2) {\scriptsize ${\uparrow}r_2$};
\draw[line width=0.5mm, cyan, ->] (45) to (56);
\draw[->,line width=0.5mm, cyan] (56) to (37) node[fill=white,inner sep=0.5pt] at (6.6,1.8) {\scriptsize ${\downarrow}r_2$};
\draw[line width=0.5mm, cyan, ->] (37) to (38);
\draw[line width=0.5mm, cyan, ->] (38) to (39);
\draw[line width=0.5mm, cyan, ->] (39) to (1);

\end{tikzpicture}
\caption{Digraph $G_{\textsc{ltc}}$ for Example \ref{example}.}
\label{fig.digraph.LTC}
\end{figure}
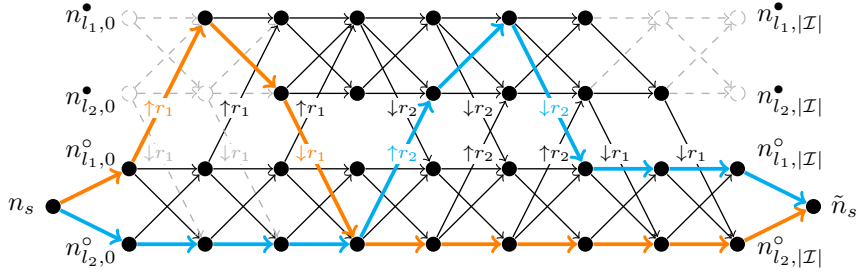

\begin{itemize}

    \item $E^{\textsc{rest}}$ has an arc $(n_{l,i}^q, n_{l,i+1}^q)$ for each $l \in L$, $i \in [|\mathcal{I}|-1]$, and $q \in \{\vacio,\lleno\}$.

    \item $E^{\textsc{trip}}$ has an arc $(n_{l_1,i}^q,n_{l_2,i + \Delta}^q)$ for each pair $(l_1,l_2) \in L \times L$ of adjacent locations, $i \in [|\mathcal{I}| - \Delta]$ with $\Delta = \VLEN(l_1,l_2)$, and $q \in \{\vacio,\lleno\}$.
    
    \item For each $r \in R$, $E^{\textsc{p},r}$ has an arc $(n_{l^\textsc{p}_r,i}^{\vacio},n_{l^\textsc{p}_r,i+s^\textsc{p}_r}^{\lleno})$ for each $i \in \mathcal{I}^\textsc{p}_r$.

    \item For each $r \in R$, $E^{\textsc{d},r}$ has an arc $(n_{l^\textsc{d}_r,i}^{\lleno},n_{l^\textsc{d}_r,i+s^\textsc{d}_r}^{\vacio})$ for each $i \in \mathcal{I}^\textsc{d}_r$.
    
    \item $E^{\textsc{sour}}$ has an arc $(n_s,n_{l,0}^{\vacio})$ for each $l \in L_V$ and $E^{\textsc{sink}}$ has an arc $(n_{l,|\mathcal{I}|}^{\vacio}, \tilde n_s)$ for each $l \in L$.

\end{itemize}
The definition of the weight vector $c$ is omitted since its adaptation is straightforward.
Fig. \ref{fig.digraph.LTC} shows the resulting construction for the instance of Example \ref{example}.
The two lower rows of nodes have the superscript ``$\vacio$'' and the two upper ones have ``$\lleno$''.
The pickup and delivery arcs are now oblique but, unlike travel arcs, they are incident on nodes with different superscripts.
For the sake of simplicity, the given definition of $G_{\textsc{ltc}}$ omitted some technicalities.
There may be nodes with the superscript ``$\lleno$'' through which it is impossible for a $(n_s,\tilde n_s)$-directed path to pass; e.g. nodes before the end time of the earliest pickup arc or after the start time of the latest delivery arc.
Such nodes and the arcs incidents on them, which are drawn with a light dashed line, can be erased from the digraph.

The greater structure of $G_{\textsc{ltc}}$ guarantees that every $(n_s,\tilde n_s)$-directed path already alternates arcs of $E^{\textsc{p}}$ and $E^{\textsc{d}}$.
Consequently, there cannot be excess capacity but services may still be unpaired or disordered.
In the worst case, $G_{\textsc{ltc}}$ has $|\mathcal{I}|.|L|$ more nodes and $|E^\textsc{rest}_{\textsc{lt}}| + |E^\textsc{trip}_{\textsc{lt}}|$ more edges than $G_{\textsc{lt}}$.
The following digraph presents an even more structured approach to modelling truck routes.

\vspace{-5pt}
\paragraph{Location--Time--Request digraph}

The weighted (simple) digraph $G_{\textsc{ltr}} \doteq (N,E,c)$ is constructed as follows.
For each $l \in L$, $i \in \mathcal{I}$, and $r \in R^{\,\vacio} \doteq R \cup \{\vacio\}$, there is a node $n_{l,i}^r$ in $N$, representing that the truck is empty in location $l$ at instant $i$ if $r=\vacio$ and is carrying request $r$ otherwise.
There are also two distinguished nodes in $N$, the source $n_s$ and the sink $\tilde n_{s}$.
The arc set $E$ follows the construction of $E_{\textsc{ltc}}$ and is the union of the following sets.

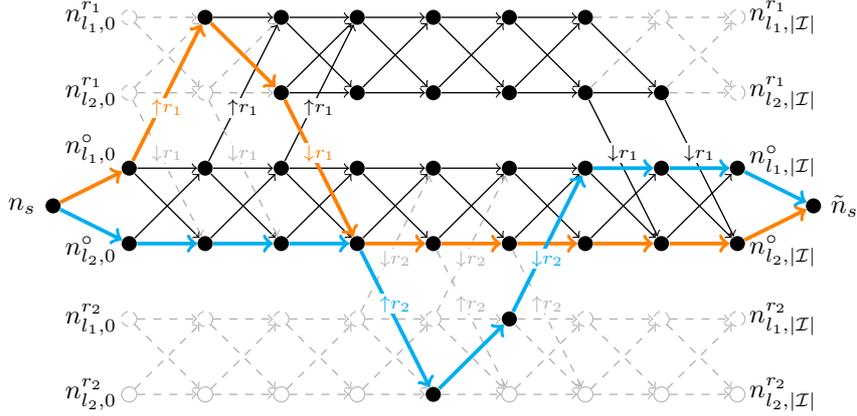
\begin{figure}
\centering
\begin{tikzpicture}

% Source and sink nodes
\node[fill,circle,inner sep=2pt, label=west:$n_s$] (0) at (0,2.5) {};
\node[fill,circle,inner sep=2pt, label=east:$\tilde n_s$] (1) at (10,2.5) {};

% Rest of the nodes
\foreach \i in {1,...,9}
{
 	\foreach \l [evaluate = {\k=\l*1}] in {2,3}
     {
 		\node[fill,circle,inner sep=2pt] (\l\i) at (\i,\k) {};
 	}       
}
\foreach \i in {1,2,3,4,6,7,8,9}
{
    \node[draw,circle,inner sep=2pt,black!30] (0\i) at (\i,0) {};
}
\node[fill,circle,inner sep=2pt] (05) at (5,0) {};
\foreach \i in {1,2,3,4,5,7,8,9}
{
 		\node[draw,circle,inner sep=2pt,black!30,densely dashed] (1\i) at (\i,1) {};
}
\node[fill,circle,inner sep=2pt] (16) at (6,1) {};
\foreach \i in {1,2,9}
{
 		\node[draw,circle,inner sep=2pt,black!30,densely dashed] (4\i) at (\i,4) {};
}
\foreach \i in {3,4,5,6,7,8}
{
 		\node[fill,circle,inner sep=2pt] (4\i) at (\i,4) {};
}
\foreach \i in {1,8,9}
{
 		\node[draw,circle,inner sep=2pt, black!30,densely dashed] (5\i) at (\i,5) {};
}
\foreach \i in {2,3,4,5,6,7}
{
 		\node[fill,circle,inner sep=2pt] (5\i) at (\i,5) {};
}

% Node labels
\node at (0.5,0) {$n^{r_2}_{l_2,0}$};
\node at (0.5,1) {$n^{r_2}_{l_1,0}$}; 
\node at (0.5,1.9) {$n^{\vacio}_{l_2,0}$};
\node at (0.5,3.25) {$n^{\vacio}_{l_1,0}$}; 
\node at (0.5,4) {$n^{r_1}_{l_2,0}$};
\node at (0.5,5) {$n^{r_1}_{l_1,0}$}; 
\node at (9.6,0) {$n^{r_2}_{l_2,|\mathcal{I}|}$};
\node at (9.6,1) {$n^{r_2}_{l_1,|\mathcal{I}|}$}; 
\node at (9.6,1.9) {$n^{\vacio}_{l_2,|\mathcal{I}|}$};
\node at (9.6,3.1) {$n^{\vacio}_{l_1,|\mathcal{I}|}$}; 
\node at (9.6,4) {$n^{r_1}_{l_2,|\mathcal{I}|}$};
\node at (9.6,5) {$n^{r_1}_{l_1,|\mathcal{I}|}$}; 

% Truck arcs
\foreach \i [evaluate={\j = int(\i+1)}] in {1,...,8}
{
    	\draw[->,line width=0.5pt] (2\i) to (3\j);
    	\draw[->,line width=0.5pt] (3\i) to (2\j);
}
\foreach \i [evaluate={\j = int(\i+1)}] in {1,2,3,4,6,7,8}
{
	\draw[->,line width=0.5pt,black!30,dashed] (0\i) to (1\j);
}
\foreach \i [evaluate={\j = int(\i+1)}] in {1,2,3,4,5,6,7,8}
{
	\draw[->,line width=0.5pt,black!30,dashed] (1\i) to (0\j);
}
\foreach \i [evaluate={\j = int(\i+1)}] in {3,4,5,6}
{
    \draw[->,line width=0.5pt] (4\i) to (5\j);
}
\foreach \i [evaluate={\j = int(\i+1)}] in {1,2,7,8}
{
	\draw[->,line width=0.5pt,black!30,dashed] (4\i) to (5\j);
}
\foreach \i [evaluate={\j = int(\i+1)}] in {3,4,5,6,7}
{
    \draw[->,line width=0.5pt] (5\i) to (4\j);
}
\foreach \i [evaluate={\j = int(\i+1)}] in {1,2,8}
{
	\draw[->,line width=0.5pt,black!30,dashed] (5\i) to (4\j);
}

% Rest arcs
\foreach \i [evaluate={\j = int(\i+1)}] in {1,...,6}
{
    \draw[->,line width=0.5pt] (3\i) to (3\j);
}
\foreach \i [evaluate={\j = int(\i+1)}] in {1,...,8}
{
    \draw[->,line width=0.5pt,black!30,dashed] (0\i) to (0\j);
    \draw[->,line width=0.5pt,black!30,dashed] (1\i) to (1\j);
}
\foreach \i [evaluate={\j = int(\i+1)}] in {1,2,8}
{
    \draw[->,line width=0.5pt,black!30,dashed] (4\i) to (4\j);
}
\foreach \i [evaluate={\j = int(\i+1)}] in {3,4,5,6,7}
{
    \draw[->,line width=0.5pt] (4\i) to (4\j);
}
\foreach \i [evaluate={\j = int(\i+1)}] in {1,7,8}
{
    \draw[->,line width=0.5pt,black!30,dashed] (5\i) to (5\j);
}
\foreach \i [evaluate={\j = int(\i+1)}] in {2,3,4,5,6}
{
    \draw[->,line width=0.5pt] (5\i) to (5\j);
}

% Pickup and delivery arcs of request 1
\foreach \i [evaluate={\j = int(\i+1)}, evaluate={\x = \i + 0.5}] in {2,3}
{
    \draw[->,line width=0.5pt] (3\i) to (5\j) node [fill=white,inner sep=0.5pt] at (\x,3.8) {\scriptsize ${\uparrow}r_1$};
}
\foreach \i [evaluate={\j = int(\i+1)}, evaluate={\x = \i + 0.5}] in {7,8}
{
    \draw[->,line width=0.5pt] (4\i) to (2\j) node[fill=white,inner sep=0.5pt] at (\x,3.2) {\scriptsize ${\downarrow}r_1$};
}
\foreach \i [evaluate={\j = int(\i+1)}, evaluate={\x = \i + 0.5}] in {1,2}
{
    \draw[->,line width=0.5pt,black!30,dashed] (4\i) to (2\j) node[fill=white,inner sep=0.5pt] at (\x,3.2) {\scriptsize ${\downarrow}r_1$};
}

% Pickup and delivery arcs of request 2
\foreach \i [evaluate={\j = int(\i+1)}, evaluate={\x = \i + 0.5}] in {5,6}
{
    \draw[->,line width=0.5pt,black!30,dashed] (2\i) to (0\j) node[fill=white,inner sep=0.5pt] at (\x,1.2) {\scriptsize ${\uparrow}r_2$};
}
\foreach \i [evaluate={\j = int(\i+1)}, evaluate={\x = \i + 0.5}] in {4,5}
{
    \draw[->,line width=0.5pt,black!30,dashed] (1\i) to (3\j) node[fill=white,inner sep=0.5pt] at (\x,1.8) {\scriptsize ${\downarrow}r_2$};
}

\draw[line width=0.5mm, orange, ->] (0) to (31);
\draw[->,line width=0.5mm, orange] (31) to (52) node[fill=white,inner sep=0.5pt] at (1.5,3.8) {\scriptsize ${\uparrow}r_1$};
\draw[line width=0.5mm, orange, ->] (52) to (43);
\draw[->,line width=0.5mm, orange] (43) to (24) node[fill=white,inner sep=0.5pt] at (3.5,3.2) {\scriptsize ${\downarrow}r_1$};
\draw[line width=0.5mm, orange, ->] (24) to (25);
\draw[line width=0.5mm, orange, ->] (25) to (26);
\draw[line width=0.5mm, orange, ->] (26) to (27);
\draw[line width=0.5mm, orange, ->] (27) to (28);
\draw[line width=0.5mm, orange, ->] (28) to (29);
\draw[line width=0.5mm, orange, ->] (29) to (1);

\draw[line width=0.5mm, cyan, ->] (0) to (21);
\draw[line width=0.5mm, cyan, ->] (21) to (22);
\draw[line width=0.5mm, cyan, ->] (22) to (23);
\draw[line width=0.5mm, cyan, ->] (23) to (24);
\draw[->,line width=0.5mm, cyan] (24) to (05) node[fill=white,inner sep=0.5pt] at (4.5,1.2) {\scriptsize ${\uparrow}r_2$};
\draw[line width=0.5mm, cyan, ->] (05) to (16);
\draw[->,line width=0.5mm, cyan] (16) to (37) node[fill=white,inner sep=0.5pt] at (6.5,1.8) {\scriptsize ${\downarrow}r_2$};
\draw[line width=0.5mm, cyan, ->] (37) to (38);
\draw[line width=0.5mm, cyan, ->] (38) to (39);
\draw[line width=0.5mm, cyan, ->] (39) to (1);

\end{tikzpicture}
\caption{Digraph $G_{\textsc{ltr}}$ for Example \ref{example}.}
\label{fig.digraph.LTR}
\end{figure}

\begin{itemize}

    \item $E^{\textsc{rest}}$ has an arc $(n_{l,i}^r, n_{l,i+1}^r)$ for each $l \in L$, $i \in [|\mathcal{I}|-1]$, and $r \in R^{\,\vacio}$.

    \item $E^{\textsc{trip}}$ has an arc $(n_{l_1,i}^r,n_{l_2,i + \Delta}^r)$ for each pair $(l_1,l_2) \in L \times L$ of adjacent locations, $i \in [|\mathcal{I}| - \Delta]$ with $\Delta = \VLEN(l_1,l_2)$, and $r \in R^{\,\vacio}$.
    
    \item For each $r \in R$, $E^{\textsc{p},r}$ has an arc $(n_{l^\textsc{p}_r,i}^{\vacio},n_{l^\textsc{p}_r,i+s^\textsc{p}_r}^r)$ for each $i \in \mathcal{I}^\textsc{p}_r$.

    \item For each $r \in R$, $E^{\textsc{d},r}$ has an arc $(n_{l^\textsc{d}_r,i}^r,n_{l^\textsc{d}_r,i+s^\textsc{d}_r}^{\vacio})$ for each $i \in \mathcal{I}^\textsc{d}_r$.
    
    \item $E^{\textsc{sour}}$ has an arc $(n_s,n_{l,0}^{\vacio})$ for each $l \in L_V$ and $E^{\textsc{sink}}$ has an arc $(n_{l,|\mathcal{I}|}^{\vacio}, \tilde n_s)$ for each $l \in L$.

\end{itemize}
For $r \in R^{\,\vacio}$, we will refer to $E^{\textsc{rest},r}$ (resp. $E^{\textsc{trip},r}$) as the subset of rest arcs (resp. trip arcs) connecting nodes with the superscript $r$.
Once again, the adaptation of the weight vector $c$ is omitted.
Fig. \ref{fig.digraph.LTR} shows the resulting construction for the instance of Example \ref{example}.
Now, the two rows of nodes in the middle have the superscript ``$\vacio$'', the upper ones have ``$r_1$'', and the lower ones have ``$r_2$''.
Again, the definition given for $G_{\textsc{ltr}}$ simplifies some technicalities in the construction.
Particularly, for all $r \in R$, the nodes with the superscript ``$r$'' before the end time of the earliest arc in $E^{\textsc{p},r}$ and the ones after the start time of the latest arc in $E^{\textsc{d},r}$ can be erased along with the arcs incidence on them.

The major benefit of $G_{\textsc{ltr}}$ is that its construction already guarantees that every $(n_s,\tilde n_s)$-directed path alternates arcs of $E^{\textsc{p}}$ and $E^{\textsc{d}}$ associated with the same request.
Thus, the services are already paired and well-ordered and there can be no excess capacity.
The main disadvantage is associated with its size, which has $|R|.|\mathcal{I}|.|L|$ more nodes and $|R|.(|E^\textsc{rest}_{\textsc{lt}}| + |E^\textsc{trip}_{\textsc{lt}}|)$ more arcs that $G_{\textsc{lt}}$ in the worst case.
The following subsection presents a digraph that can be utilised to represent the driver routes.

\subsection{Driver routes} \label{sec.models.driver.routes}

\vspace{-5pt}
\paragraph{Location--Time digraph with taxi arcs} 

Let $G_{\textsc{ltx}} \doteq (N,E,c)$ be the weighted multidigraph that has the same node set as $G_{\textsc{lt}}$, i.e. $N = N_{\textsc{lt}}$, and whose arc set is the union of the following sets.
\begin{itemize}
    \item $E^{\textsc{rest}}$, $E^{\textsc{trip}}$, $E^{\textsc{p}}$, $E^{\textsc{r}}$, and $E^{\textsc{sink}}$, as defined in $G_{\textsc{lt}}$.
    \item $E^{\textsc{taxi}}$ has an arc $(n_{l_1,i},n_{l_2,i + \Delta})$ for each pair $(l_1,l_2) \in L \times L$ of adjacent locations and $i \in [|\mathcal{I}| - \Delta]$ with  $\Delta = \TLEN(l_1,l_2)$, representing that a driver takes a taxi.
    \item $E^{\textsc{sour}}$ has an arc $(n_s,n_{l,0}^{\vacio})$ for each $l \in L_D$.
\end{itemize}
The weight of an arc $e = (n_{l_1,i},n_{l_2,i + \Delta}) \in E^{\textsc{taxi}}$ is the taxi travel cost between $l_1$ and $l_2$, i.e. $c(e) = \TCOST(l_1,l_2)$, and the remaining are zero-weighted arcs.
Fig. \ref{fig.digraph.LTX} shows the resulting construction for the instance of Example \ref{example} and the driver routes of Fig. \ref{fig.example}(a).
Since this example considers travel times to be identical for trucks and taxis, the arcs in $E^{\textsc{taxi}}$ are parallel to those in $E^{\textsc{trip}}$.
To distinguish them, taxi arcs are drawn curved with the label ``\textsc{t}''.
The green $(n_s,\tilde n_s)$-directed path corresponds to the route of $d_1$ and the pink one to $d_2$.

\begin{figure}
\centering
\begin{tikzpicture}

% Source and sink nodes
\node[fill,circle,inner sep=2pt, label=west:$n_s$] (0) at (0,0.5) {};
\node[fill,circle,inner sep=2pt, label=east:$\tilde n_s$] (1) at (12.5,0.5) {};

% Rest of the nodes
\foreach \i [evaluate={\x = 1.25*\i}] in {1,...,9}
{
 	\foreach \l [evaluate={\y = \l/1.5}] in {0,2}
     {
 		\node[fill,circle,inner sep=2pt] (\l\i) at (\x,\y) {};
 	}       
}

% Node labels
\node at (0.5,-0.2) {$n_{l_2,0}$};
\node at (0.5,1.4) {$n_{l_1,0}$}; 

% Source and sink arcs
\draw[->,line width=0.5pt] (09) to (1);  

% Truck arcs
\foreach \i [evaluate={\j = int(\i+1)}] in {1,2,3,4,6,7,8}
{
    \draw[->,line width=0.5pt] (0\i) to (2\j);
}
\foreach \i [evaluate={\j = int(\i+1)}] in {1,3,4,5,6,7,8}
{
    \draw[->,line width=0.5pt] (2\i) to (0\j);
}

% Taxi arcs
\foreach \i [evaluate={\j = int(\i+1)}, evaluate={\x = 1.25*\i + 1.15}] in {1,2,3,5,6,7,8}
{
    \draw[->,line width=0.5pt, bend right=30] (0\i) to (2\j) node[fill=white,inner sep=1pt] at (\x,0.95) {\footnotesize \textsc{t}};
}
\foreach \i [evaluate={\j = int(\i+1)}, evaluate={\x = 1.25*\i + 1.15}] in {1,...,8}
{
    \draw[->,line width=0.5pt, bend left=30] (2\i) to (0\j) node[fill=white,inner sep=1pt] at (\x,0.4) {\footnotesize \textsc{t}};  
}

% Rest arcs
\foreach \i [evaluate={\j = int(\i+1)}] in {4,5,...,8}
{
    \draw[->,line width=0.5pt] (0\i) to (0\j);
}
\foreach \i [evaluate={\j = int(\i+1)}] in {1,...,4}
{
    \draw[->,line width=0.5pt] (2\i) to (2\j);
}

% Pickup and delivery arcs of request 1
\foreach \i [evaluate={\j = int(\i+1)}, evaluate={\x = 1.25*\i + 0.65}] in {2,3}
{
    \draw[->,line width=0.5pt, bend left = 90, looseness = 1.25] (2\i) to (2\j) node [fill=white, inner sep=1pt] at (\x,1.9) {\footnotesize ${\uparrow}r_1$};
}
\foreach \i [evaluate={\j = int(\i+1)}, evaluate={\x = 1.25*\i + 0.65}] in {1,2,7,8}
{
    \draw[->,line width=0.5pt, bend right = 90, looseness = 1.25] (0\i) to (0\j) node [fill=white, inner sep=1pt] at (\x,-0.5) {\footnotesize ${\downarrow}r_1$};
}

% Pickup and delivery arcs of request 2
\foreach \i [evaluate={\j = int(\i+1)}, evaluate={\x = 1.25*\i + 0.65}] in {5,6}
{
    \draw[->,line width=0.5pt, bend right = 90, looseness = 1.25] (0\i) to (0\j) node [fill=white, inner sep=1pt] at (\x,-0.5) {\footnotesize ${\uparrow}r_2$};
}
\foreach \i [evaluate={\j = int(\i+1)}, evaluate={\x = 1.25*\i + 0.65}] in {4,5}
{
    \draw[->,line width=0.5pt, bend left = 90, looseness = 1.25] (2\i) to (2\j) node [fill=white, inner sep=1pt] at (\x,1.9) {\footnotesize ${\downarrow}r_2$};
}

\draw[line width=0.5mm, ForestGreen, ->] (0) to (21);
\draw[->,line width=0.5mm, bend left = 90, looseness = 1.25, ForestGreen] (21) to (22) node [fill=white, inner sep=1pt] at (1.9,1.9) {\footnotesize ${\uparrow}r_1$};
\draw[line width=0.5mm, ForestGreen, ->] (22) to (03);
\draw[line width=0.5mm, ForestGreen, ->] (03) to (04);
\draw[->,line width=0.5mm, bend right = 90, looseness = 1.25, ForestGreen] (04) to (05) node [fill=white, inner sep=1pt] at (5.65,-0.5) {\footnotesize ${\uparrow}r_2$};
\draw[line width=0.5mm, ForestGreen, ->] (05) to (26);
\draw[line width=0.5mm, ForestGreen, ->] (26) to (27);
\draw[line width=0.5mm, ForestGreen, ->] (27) to (28);
\draw[line width=0.5mm, ForestGreen, ->] (28) to (29);
\draw[line width=0.5mm, ForestGreen, ->] (29) to (1);

\draw[line width=0.5mm, magenta, ->] (0) to (01);
\draw[line width=0.5mm, magenta, ->] (01) to (02);
\draw[line width=0.5mm, magenta, ->] (02) to (03);
\draw[->,line width=0.5mm, bend right = 90, looseness = 1.25, magenta] (03) to (04) node [fill=white, inner sep=1pt] at (4.4,-0.5) {\footnotesize ${\downarrow}r_1$};
\draw[line width=0.5mm, magenta, ->, bend right=30] (04) to (25) node[fill=white,inner sep=1pt] at (6.15,0.95) {\footnotesize \textsc{t}};
\draw[line width=0.5mm, magenta, ->] (25) to (26);
\draw[->,line width=0.5mm, bend left = 90, looseness = 1.25, magenta] (26) to (27) node [fill=white, inner sep=1pt] at (8.15,1.9) {\footnotesize ${\downarrow}r_2$};
\draw[line width=0.5mm, magenta, dashed] (27) to (28);
\draw[line width=0.5mm, magenta, dashed] (28) to (29);
\draw[line width=0.5mm, magenta, dashed] (29) to (1);

\end{tikzpicture}
\caption{Digraph $G_{\textsc{ltx}}$ for Example \ref{example}.}
\label{fig.digraph.LTX}
\end{figure}

The construction of $G_{\textsc{ltx}}$ guarantees that every driver route corresponds to a $(n_s,\tilde n_s)$-directed path, but it is also necessary to forbid some paths for the converse to hold.
Let $P$ be a $(n_s,\tilde n_s)$-directed path, there are two categories of forbidden paths according to the labour laws of the case study. 

\begin{itemize}
\item \textit{Daily rest.} Since drivers must rest at least 12 hours in each 24-hour interval, then for all $i \in [I(H-1)]$ and $P'$ subdirected path of $P$ from nodes $n_{l_1,i}$ to $n_{l_2,i + I}$ with $l_1,l_2 \in L$, the total number of arcs of $E^{\textsc{rest}}$ in $P'$ cannot be lower than $\frac{I}{2}$. 
\item \textit{Weekly rest.} Since drivers must have at least 1 day off in each 7-day interval, then for all $j \in [H-7]$ and $P'$ subdirected path of $P$ from nodes $n_{l_1,Ij}$ to $n_{l_2,I(j+7)}$ with $l_1,l_2 \in L$, there must be a subdirected path of $P'$ with all its arcs in $E^{\textsc{rest}}$ from nodes $n_{l,Ij'}$ to $n_{l,I(j'+1)}$ for some $l \in L$ and $j' \in \{j,\ldots,j+6\}$.
\end{itemize}

Unlike the approach followed for trucks, no further representations are proposed for driver routes, since enriching the nodes with cargo information has no clear advantage for them. 
In what follows, we will discuss how to synchronise the routes of the trucks and drivers.

\subsection{Synchronisation of routes} \label{sec.models.sync.routes}

Synchronization is of utmost importance to ensure that trucks and drivers match in time and space. 
The selected representation for the routes has a significant benefit as it enables synchronization to be expressed in a simplified manner.
The most direct case is when truck routes are represented with $G_{\textsc{lt}}$ and driver routes with $G_{\textsc{ltx}}$, since the arc set of the former is contained in the latter.
For each truck route that passes through a trip, pickup, or delivery arc (which are the ones that require synchronization), there must be some driver route that passes through that same arc, and vice-versa.
However, since drivers can travel as passengers, at most two of them can share the same arc of a truck route.
Formally, for each $e \in E^{\textsc{trip}}_{\textsc{ltx}} \cup E^{\textsc{p}}_{\textsc{ltx}} \cup E^{\textsc{d}}_{\textsc{ltx}}$, if $V_e$ and $D_e$ are the numbers of trucks and drivers whose routes pass through $e$, respectively, then it must hold $V_e \leq D_e \leq 2V_e.$

For the other representations, the arc correspondence between the digraphs is less direct and no longer one-to-one.
For example, a driver can travel in either an empty or a full truck, but $G_{\textsc{ltc}}$ and $G_{\textsc{ltr}}$ identify such trip with different arcs (varying the superscript of the nodes they connect).
To overcome this, functions $\arcLTC$ and $\arcLTR$ are defined, that return the subset of arcs of $E_{\textsc{ltc}}$ and $E_{\textsc{ltr}}$ that match a given $e \in E_{\textsc{ltx}}$, respectively.
For the sake of brevity, we present only the relevant cases of $\arcLTC$, even though the definition can be extended to cover its entire domain, as well as extended to $\arcLTR$.
\begin{eqnarray*}
\arcLTC(e) \doteq \span
\begin{cases}
%\{e' \in E^{\textsc{rest}}_{\textsc{ltc}} : e' = (n^q_{l,i_1},n^q_{l,i_2}) \wedge q \in \{\vacio,\lleno\}\} & \text{if } e = (n_{l,i_1},n_{l,i_2}) \in E^{\textsc{rest}}_{\textsc{ltx}}\\ 
\{e' \in E^{\textsc{trip}}_{\textsc{ltc}} : e' = (n^q_{l_1,i_1},n^q_{l_2,i_2}) \wedge q \in \{ \vacio,\lleno\}\} & \text{if } e = (n_{l_1,i_1},n_{l_2,i_2}) \in E^{\textsc{trip}}_{\textsc{ltx}}.\\ 
\{e' \in E^{\textsc{p},r}_{\textsc{ltc}} : e' = (n^{\vacio}_{l,i_1},n^{\lleno}_{l,i_2})\} & \text{if } e = (n_{l,i_1},n_{l,i_2}) \in E^{\textsc{p},r}_{\textsc{ltx}}.\\
\{e' \in E^{\textsc{d},r}_{\textsc{ltc}}: e' = (n^{\lleno}_{l,i_1},n^{\vacio}_{l,i_2})\} & \text{if } e = (n_{l,i_1},n_{l,i_2}) \in E^{\textsc{d},r}_{\textsc{ltx}}.\\
%\{e' \in E^{\textsc{sour}}_{\textsc{ltc}} : e' = (n_f, n_{l,0}^{\vacio})\} & \text{if } e = (n_f, n_{l,0}) \in E^{\textsc{sour}}_{\textsc{ltx}}\\
%\{e' \in E^{\textsc{sink}}_{\textsc{ltc}} : e' = (n_{l,I.H}^{\vacio},n_s)\} & \text{if } e = (n_{l,I.H},n_s) \in E^{\textsc{sink}}_{\textsc{ltx}}\\
%\{\} & \text{if } e \in E^{\textsc{taxi}}_{\textsc{ltx}}
\end{cases}
\end{eqnarray*}  
Then, the previous inequality can be stated as
$ \sum_{e' \in \arcLTC(e)} V_{e'}  \leq D_e \leq 2 \sum_{e' \in \arcLTC(e)} V_{e'}.$

\section{Integer linear programs}\label{sec.ilp}

Before presenting the ILP formulations, some additional notations must be introduced. 
Given a node $n \in N$, $\Gamma^-(n) \subset E$ (resp. $\Gamma^+(n)$) is the subset of incoming (resp. outgoing) arcs of $n$. 
Given $n_{l,i} \in N$ and $e \in \Gamma^+(n_{l,i})$ (resp. $e \in \Gamma^-(n_{l,i})$), then $\origen(e) \doteq l$ and $\inicio(e) \doteq i$ (resp. $\destino(e) \doteq l$ and $\fin(e) \doteq i$).
Given $v \in V$ and $d \in D$, $E^v \subset E$ and $E^d \subset E$ are the subsets of possible arcs for $v$ and $d$ according to their start location, resp., i.e. $E^v \doteq E \setminus \{e \in E^{\textsc{sour}} : \destino(e) \neq l_v \}$ and $E^d \doteq E \setminus \{e \in E^{\textsc{sour}} : \destino(e) \neq l_d \}$.

\vspace{-5pt}
\paragraph{LT formulation} 
The first ILP formulation is based on $G_{\textsc{lt}}$ and $G_{\textsc{ltx}}$ to represent truck and driver routes, respectively.
For each $v \in V$ and $e \in E^v_{\textsc{lt}}$, there is a binary variable $X_{ve}$ such that $X_{ve} = 1$ if and only if $v$ passes through $e$.
For each $d \in D$ and $e \in E^d_{\textsc{ltx}}$, there is a binary variable $Y_{de}$ such that $Y_{de} = 1$ if and only if $d$ passes through $e$.
For each $d \in D$ and $j \in [H-1]$, there is a binary variable $W_{dj}$ such that $W_{dj} = 1$ if $d$ has day $j$ off.
{
\allowdisplaybreaks
\begin{align}
\text{(LT)} \hspace{15pt} & \min \sum_{v \in V} \sum_{e \in E^v_{\textsc{lt}}} c_{\textsc{lt}}(e) X_{ve} + \sum_{d \in D} \sum_{e \in E^d_{\textsc{ltx}}} c_{\textsc{ltx}}(e) Y_{de} \span\span \label{form.LT.LTX.Func.Obj}  \\[-2pt] 
\text{s.t.} \hspace{20pt} & \sum_{\mathclap{e \in \Gamma^-(n) \cap E^v_{\textsc{lt}}}} \hspace{7pt} X_{ve} = \hspace{7pt} \sum_{\mathclap{e \in \Gamma^+(n) \cap E^v_{\textsc{lt}}}} \hspace{5pt} X_{ve} & \forall v \in V,\ n \in N_{\textsc{lt}} \setminus \{n_s,\tilde n_s\}, \label{form.LT.LTX.Truck.Flow} \\[-2pt] 
& \sum_{\mathclap{e \in \Gamma^-(n) \cap E^d_{\textsc{ltx}}}} \hspace{7pt} Y_{de} = \hspace{7pt} \sum_{\mathclap{e \in \Gamma^+(n) \cap E^d_{\textsc{ltx}}}} \hspace{5pt} Y_{de} & \forall d \in D,\ n \in N_{\textsc{ltx}} \setminus \{n_s,\tilde n_s\}, \label{form.LT.LTX.Driver.Flow} \\[-2pt] 
& \sum_{v \in V} \sum_{e \in E^{\textsc{p},r}_{\textsc{lt}}} X_{ve} = 1 & \forall r \in R, \label{form.LT.LTX.Truck.Repeated} \\[-2pt] 
& \sum_{\mathclap{e \in E^{\textsc{p},r}_{\textsc{lt}}}} X_{ve} = \sum_{e \in E^{\textsc{d},r}_{\textsc{lt}}} X_{ve} & \forall v \in V,\ r \in R, \label{form.LT.LTX.Truck.Unpaired} \\[-2pt]
& \sum_{v \in V} \hspace{8pt} \sum_{\mathclap{\substack{e' \in E^{\textsc{p},r}_{\textsc{lt}}:\\ \fin(e') \leq \inicio(e)}}} \hspace{3pt} X_{ve'} \geq \sum_{v \in V} X_{ve} & \forall r \in R,\ e \in E^{\textsc{d},r}_{\textsc{lt}}, \label{form.LT.LTX.Truck.Precedence} \\[-2pt]
& \sum_{r \in R} \hspace{4pt} \Bigl( \sum_{\mathclap{\substack{e \in E^{\textsc{p},r}_{\textsc{lt}}:\\ \fin(e) \leq i}}} X_{ve} - \sum_{\mathclap{\substack{e \in E^{\textsc{d},r}_{\textsc{lt}}:\\ \fin(e) \leq i}}} X_{ve} \Bigr) \leq 1 & \forall v \in V,\ i \in {\textstyle \bigcup_{r \in R} \mathcal{I}^\textsc{p}_r}, \label{form.LT.LTX.Truck.Capacity} \\[-2pt]
& \sum_{\mathclap{\substack{e \in E^{\textsc{rest}}_{\textsc{ltx}}:\\i \leq \inicio(e) < i+I}}} \hspace{5pt} Y_{de} \geq I/2 & \forall d \in D,\ i \in  [I(H-1)], \label{form.LT.LTX.Driver.Daily} \\[-2pt]
& \sum_{\mathclap{j \leq j' \leq j+6}} \hspace{5pt} W_{dj'} \geq 1 & \forall d \in D,\ j \in [H-7], \label{form.LT.LTX.Driver.Weekly.I} \\[-2pt]
& \sum_{\mathclap{\substack{e \in E^{\textsc{desc}}_{\textsc{ltx}}:\\Ij \leq \inicio(e) < I(j+1)}}} \hspace{5pt} Y_{de} \geq IW_{dj} & \forall d \in D,\ j \in [H-1], \label{form.LT.LTX.Driver.Weekly.II} \\[-2pt]
& \sum_{v \in V} X_{ve} \leq \sum_{d \in D} Y_{de} \leq 2\sum_{v \in V} X_{ve} & \forall e \in E^{\textsc{trip}}_{\textsc{ltx}} \cup E^{\textsc{p}}_{\textsc{ltx}} \cup E^{\textsc{d}}_{\textsc{ltx}}, \label{form.LT.LTX.Synchronisation}  \\[-2pt]
& X_{ve} \in \{0,1\} & \forall v \in V,\ e \in E^v_{\textsc{lt}}, \label{form.LT.LTX.Range.I} \\[-2pt]
& Y_{de} \in \{0,1\} & \forall d \in D,\ e \in E^c_{\textsc{ltx}}, \label{form.LT.LTX.Range.II} \\[-2pt]
& W_{dj} \in \{0,1\} & \forall d \in D,\ j \in [H-1]. \label{form.LT.LTX.Range.III}
\end{align}
}
The objective function (\ref{form.LT.LTX.Func.Obj}) minimises total travel costs and penalties for late deliveries.
Constraints (\ref{form.LT.LTX.Truck.Flow})--(\ref{form.LT.LTX.Driver.Flow}) are the flow conservation equations for trucks and drivers, respectively, and they guarantee that each route is empty or a $(n_s, \tilde n_s)$-directed path. 
Then, (\ref{form.LT.LTX.Truck.Repeated})--(\ref{form.LT.LTX.Truck.Capacity}) and (\ref{form.LT.LTX.Driver.Daily})--(\ref{form.LT.LTX.Driver.Weekly.II}) prohibit directed paths that do not correspond to truck and driver routes according to the categories of Section \ref{sec.models.truck.routes} and \ref{sec.models.driver.routes}, respectively.
Observe that (\ref{form.LT.LTX.Truck.Repeated}) also prevent a request from being executed by more than one truck and adding similar constraints for the deliveries is not necessary since they are implied.
%by (\ref{form.LT.LTX.Truck.Repeated}) and (\ref{form.LT.LTX.Truck.Unpaired}).
When considering individual truck routes, 
(\ref{form.LT.LTX.Truck.Repeated}) and (\ref{form.LT.LTX.Truck.Precedence}) disaggregated by truck are also implied.
While (\ref{form.LT.LTX.Driver.Weekly.I}) force each driver to have at least one day off, (\ref{form.LT.LTX.Driver.Weekly.II}) guarantee that their route contains only rest arcs on such day.
Constraints (\ref{form.LT.LTX.Synchronisation}) synchronise the routes.
Finally, (\ref{form.LT.LTX.Range.I})--(\ref{form.LT.LTX.Range.III}) are the integrality constraints. 

\vspace{-5pt}
\paragraph{LTC formulation} 
In this ILP formulation, $G_{\textsc{ltc}}$ is used to represent the truck routes.
In addition to the above variables $Y$ and $W$, for each $v \in V$ and $e \in E^v_{\textsc{ltc}}$, there is a binary variable $X_{ve}$ such that $X_{ve} = 1$ if and only if $v$ passes through $e$.
{
\allowdisplaybreaks
\begin{align}
\text{(LTC)} \hspace{15pt} & \min \sum_{v \in V} \sum_{e \in E^v_{\textsc{ltc}}} c_{\textsc{ltc}}(e) X_{ve} + \sum_{d \in D} \sum_{e \in E^d_{\textsc{ltx}}} c_{\textsc{ltx}}(e) Y_{de} \span\span \label{form.LTC.LTX.Func.Obj}  \\[-2pt]  
\text{s.t.} \hspace{20pt} &
(\ref{form.LT.LTX.Driver.Flow}), (\ref{form.LT.LTX.Driver.Daily}), (\ref{form.LT.LTX.Driver.Weekly.I}), (\ref{form.LT.LTX.Driver.Weekly.II}), (\ref{form.LT.LTX.Range.II}), (\ref{form.LT.LTX.Range.III}), \span\span\nonumber \\[-2pt] 
& \sum_{\mathclap{e \in \Gamma^-(n) \cap E^v_{\textsc{ltc}}}} \hspace{7pt} X_{ve} = \hspace{7pt} \sum_{\mathclap{e \in \Gamma^+(n) \cap E^v_{\textsc{ltc}}}} \hspace{5pt} X_{ve} & \forall v \in V,\ n \in N_{\textsc{ltc}} \setminus \{n_s,\tilde n_s\}, \label{form.LTC.LTX.Truck.Flow} \\[-2pt] 
& \sum_{v \in V} \sum_{e \in E^{\textsc{p},r}_{\textsc{ltc}}} X_{ve} = 1 & \forall r \in R, \label{form.LTC.LTX.Truck.Repeated} \\[-2pt] 
& \sum_{\mathclap{e \in E^{\textsc{p},r}_{\textsc{ltc}}}} X_{ve} = \sum_{e \in E^{\textsc{d},r}_{\textsc{ltc}}} X_{ve} & \forall v \in V,\ r \in R, \label{form.LTC.LTX.Truck.Unpaired} \\[-2pt] 
& \sum_{v \in V} \hspace{8pt} \sum_{\mathclap{\substack{e' \in E^{\textsc{p},r}_{\textsc{ltc}}:\\ \fin(e') \leq \inicio(e)}}} \hspace{3pt} X_{ve'} \geq \sum_{v \in V} X_{ve} & \forall r \in R,\ e \in E^{\textsc{d},r}_{\textsc{ltc}}, \label{form.LTC.LTX.Truck.Precedence} \\[-2pt] 
& \sum_{v \in V} \hspace{3pt} \sum_{\mathclap{\substack{e' \in\\ \arcLTC(e)}}} X_{ve'} \leq \sum_{d \in D} Y_{de} \leq 2\sum_{v \in V} \hspace{3pt} \sum_{\mathclap{\substack{e' \in\\ \arcLTC(e)}}} X_{ve'} & \forall e \in E^{\textsc{trip}}_{\textsc{ltx}} \cup E^{\textsc{p}}_{\textsc{ltx}} \cup E^{\textsc{d}}_{\textsc{ltx}}, \label{form.LTC.LTX.Synchronisation}  \\[-2pt] 
& X_{ve} \in \{0,1\} & \forall v \in V,\ e \in E^v_{\textsc{ltc}}. \label{form.LTC.LTX.Range.I}
\end{align}
}

\vspace{-5pt}
\paragraph{LTR formulation} 
Finally, the ILP formulation is presented that considers $G_{\textsc{ltr}}$ to represent the truck routes.
In addition to the above variables $Y$ and $W$, for each $e \in E_{\textsc{ltr}}$, there is an integer variable $X_{e}$ whose value is the number of trucks that passes through $e$.
The range of $X_e$ is upper bounded by the capacity of $e$.
Particularly, the number of trucks in each start location defines the capacity for each outgoing arc of the source node, the capacity for the arcs where the truck is loaded is one to avoid repeated services, otherwise the capacity is limited by the total number of trucks.
{
\allowdisplaybreaks
\begin{align}
\text{(LTR)} \hspace{15pt} & \min \sum_{e \in E_{\textsc{ltr}}} c_{\textsc{ltr}}(e) X_{e} + \sum_{d \in D} \sum_{e \in E^d_{\textsc{ltx}}} c_{\textsc{ltx}}(e) Y_{de} \span \label{form.LTR.LTX.Func.Obj}  \\[-2pt]  
\text{s.t.} \hspace{20pt} & (\ref{form.LT.LTX.Driver.Flow}), (\ref{form.LT.LTX.Driver.Daily}), (\ref{form.LT.LTX.Driver.Weekly.I}), (\ref{form.LT.LTX.Driver.Weekly.II}), (\ref{form.LT.LTX.Range.II}),(\ref{form.LT.LTX.Range.III}),\span\nonumber \\[-2pt] 
& \sum_{\mathclap{e \in \Gamma^-(n)}} X_{e} = \sum_{\mathclap{e \in \Gamma^+(n)}} X_{e} & \forall n \in N_{\textsc{ltr}} \setminus \{n_s,\tilde n_s\}, \label{form.LTR.LTX.Truck.Flow} \\[-2pt] 
& \sum_{\mathclap{e \in E^{\textsc{p},r}_{\textsc{ltr}}}} X_{e} = 1 & \forall r \in R, \label{form.LTR.LTX.Truck.Repeated} \\[-2pt] 
& \sum_{\mathclap{\substack{e' \in \arcLTR(e)}}} X_{e'} \leq \sum_{d \in D} Y_{de} \leq 2\sum_{\mathclap{\substack{e' \in \arcLTR(e)}}} X_{e'} & \forall e \in E^{\textsc{trip}}_{\textsc{ltx}} \cup E^{\textsc{p}}_{\textsc{ltx}} \cup E^{\textsc{d}}_{\textsc{ltx}}, \label{form.LTR.LTX.Synchronisation}  \\[-2pt] 
& X_{e} \in \{0,\ldots,\capacity(e)\} & \forall e \in E_{\textsc{ltr}}. \label{form.LTR.LTX.Range.I}
\end{align}
}
Where the function $\capacity$ is defined as $\capacity(e) \doteq |\{v \in V: l_v = l\}|$ if $e \in E^{\textsc{sour}}_{\textsc{ltr}} \text{ and } \destino(e) = l$, $\capacity(e) \doteq 1$ if, for some $r \in R, e \in E^{\textsc{p},r}_{\textsc{ltr}} \cup E^{\textsc{d},r}_{\textsc{ltr}} \cup E^{\textsc{rest},r}_{\textsc{ltr}} \cup E^{\textsc{trip},r}_{\textsc{ltr}}$, and $\capacity(e) \doteq |V|$ otherwise.
%\begin{eqnarray*} 
%\capacity(e) \doteq
%\begin{cases}
%|\{v \in V: l_v = l\}| & \text{if } e \in E^{\textsc{sour}}_{\textsc{ltr}} \text{ and } \destino(e) = l. \\
%1 & \text{if, for some } r \in R, e \in E^{\textsc{p},r}_{\textsc{ltr}} \cup E^{\textsc{d},r}_{\textsc{ltr}} \cup E^{\textsc{rest},r}_{\textsc{ltr}} \cup E^{\textsc{trip},r}_{\textsc{ltr}}.\\
%|V| & \text{otherwise}. 
%\end{cases}
%\end{eqnarray*}
In this case, translating a feasible solution to truck routes is not straightforward. 
It involves partitioning the flow given by $X$ into $(n_s,\tilde n_s)$-directed paths, which can be achieved in polynomial time.

\section{Alternative constraints and valid inequalities}\label{sec.valid.ineq}

This section describes some alternative constraints and valid inequalities that cut fractional solutions, which will later be used as cutting planes.
For simplicity, they are expressed for the LT formulation, but they can be adapted to the others, except in specific cases that will be commented on in due course.

\vspace{-5pt}
\paragraph{Precedence inequalities}
They are inspired by an alternative definition of the forbidden paths mentioned in Section \ref{sec.models.truck.routes}. Let $P$ be a $(n_s,\tilde n_s)$-directed path,
\begin{itemize}
\item \textit{Disordered services (alternative).} For all $r \in R$ and $P'$ subdirected path of $P$ from $n_s$, the difference between the total number of arcs of $E^{\textsc{d},r}$ and $E^{\textsc{p},r}$ in $P'$ cannot be greater than zero.
\end{itemize}
This alternative definition now excludes, for example, a path that has an arc of $E^{\textsc{p},r}$ followed by two arcs of  $E^{\textsc{d},r}$.
Despite being more restrictive than the original, they coincide when $P$ has no repeated services.
%However, such paths also have repeated services ($r$ is delivered twice), so they are forbidden anyway.
At the same time, it is unnecessary to consider all arcs of $E^{\textsc{p},r}$ in $P'$, but only those that give the truck enough time to trip to the delivery location.
Then, the following valid inequalities arise:
\begin{align}
& \hspace{10pt} \sum_{\mathclap{\substack{e \in E^{\textsc{p},r}_{\textsc{lt}}:\\[2pt] \fin(e) \leq i - \TLEN(l^{\textsc{p}}_r, l^{\textsc{d}}_r)}}} \hspace{10pt} \sum_{v \in V} X_{ve} \geq \sum_{\mathclap{\substack{e \in E^{\textsc{d},r}_{\textsc{lt}}:\\[2pt] \inicio(e) \leq i}}} \hspace{10pt} \sum_{v \in V} X_{ve} & \forall r \in R,\ i \in \mathcal{I}^{\textsc{d}}_r. \tag{PREC}
\end{align}
An additional advantage is that they can replace constraints (\ref{form.LT.LTX.Truck.Precedence}), yielding a tighter ILP formulation.

\vspace{-5pt}
\paragraph{Synchronisation constraints} 
The synchronisation of routes requires one or two drivers to be present during truck activities, as proposed in Section \ref{sec.models.sync.routes}.
While it may seem logical for drivers to avoid taking taxis by travelling as passengers, there is no advantage to sharing a truck in other circumstances.
Observe that during loading or unloading, one of the drivers may always be forced to remain in the truck and the other to rest outside, which preserves the cost and feasibility of the solution.
The following constraints can achieve this restriction on the solution space:
\begin{align*}
\sum_{v \in V} X_{ve} & = \sum_{d \in D} Y_{de} & \forall e \in E^{\textsc{p}}_{\textsc{ltx}} \cup E^{\textsc{d}}_{\textsc{ltx}}. \tag{SYNC$_1$}
\end{align*}
The opposite approach involves allowing an unlimited number of drivers in the truck during loading and unloading.
In this case, this relaxation is achieved by the following constraints:
\begin{align*}
\sum_{v \in V} X_{ve} & \leq \sum_{d \in D} Y_{de} & \forall e \in E^{\textsc{p}}_{\textsc{ltx}} \cup E^{\textsc{d}}_{\textsc{ltx}}. \tag{SYNC$_2$}
\end{align*}
The constraints (\ref{form.LT.LTX.Synchronisation}) on pickup and delivery arcs can be replaced by any of the above.
In particular, SYNC$_1$ eliminate some feasible solutions but keep some optimal ones, and SYNC$_2$ introduce additional feasible solutions but an optimal one can always be recovered for the original problem.

\vspace{-5pt}
\paragraph{Pickup-and-delivery trip inequalities}
Their motivation comes from observing that trucks always travel from the pickup location to the delivery location (not necessarily directly) for each request they serve.
The first family of inequalities forces a truck route to have as many trip arcs with origin (resp. destination) in a given location as requests it loads (resp. unloads) there.
\begin{align}
\sum_{\substack{e \in E^\textsc{trip}_{\textsc{lt}}:\\ \origen(e) = l}} X_{ve} & \geq \sum_{\substack{r \in R:\\ l^{\textsc{p}}_r = l}} \hspace{2pt} \sum_{e \in E^{\textsc{p},r}_{\textsc{lt}}} X_{ve} & \forall l \in L,\ v \in V. \tag{PD$_1$}
\end{align}
To avoid making the description too long, similar inequalities for the deliveries are also included. %despite not being displayed.
However, it is important to note that such trips always occur after loading (resp. before unloading).
By restricting the arcs according to their start time, the following alternative inequalities are obtained.
\begin{align}
\hspace{5pt} \sum_{\mathclap{\substack{e \in E^\textsc{trip}_{\textsc{lt}}:\\ \origen(e) = l,\ \inicio(e) \geq i+s^{\textsc{p}}_r}}} \hspace{5pt} X_{ve} & \geq \sum_{\substack{r \in R:\\ l^{\textsc{p}}_r = l}} \hspace{2pt} \sum_{\substack{e \in E^{\textsc{p},r}_{\textsc{lt}}:\\ \inicio(e) \geq i}} X_{ve} & \forall l \in L,\ i  \in \bigcup_{\substack{r \in R:\\ l^{\textsc{p}}_r = l}} \mathcal{I}^{\textsc{p}}_r,\ v \in V. \tag{PD$_2$}
\end{align}
Observe that PD$_2$ include PD$_1$ when $i$ is the earliest time instant.
%A different approach consists of simultaneously limiting the start time of the trip both below and above.
The last family of these inequalities requires the truck to travel from the pickup location (resp. to the delivery location) both after loading and before unloading a given request.
Instead of considering a single truck, the following inequalities are generalised over a subset $V'$ of trucks.
\begin{align*}
\hspace{-20pt} \sum_{v \in V'} \hspace{5pt} \Bigl( \sum_{\mathclap{\substack{e \in E^{\textsc{p},r}_{\textsc{lt}}:\\ \inicio(e_1) \leq \inicio(e)}}} X_{ve} +  \sum_{\mathclap{\substack{e \in E^{\textsc{d},r}_{\textsc{lt}}:\\ \fin(e) \leq \fin(e_2)}}} X_{ve} \Bigr) - 1 & \leq \sum_{v \in V'} \hspace{7pt} \sum_{\mathclap{\substack{e \in E^{\textsc{trip}}_{\textsc{lt}}:\\ \origen(e) = l^\textsc{p}_r,\\ \fin(e_1) \leq \inicio(e),\ \fin(e) \leq \inicio(e_2)}}} X_{ve} & \forall~ r \in R,\ e_1 \in E^{\textsc{p},r},\ e_2 \in E^{\textsc{d},r},\ V' \subset V. \tag{PD$_3$}
\end{align*}
They should be taken into consideration only when $\fin(e_1) + \TLEN(l^{\textsc{r}}_p, l^{\textsc{e}}_p) \leq \inicio(e_2)$.
The validity follows from the fact that $r$ can be served by at most one truck in $V'$ since repeated services are prohibited, thus one is an upper bound for the left-hand side.

In their adaptations to the LTC formulation, a significant optimization can be applied. 
Since the trip from the pickup location (resp. to the delivery location) must occur while the truck is loaded, the set $E^{\textsc{trip}}_{\textsc{ltc}}$ in the summation can be tightened to $E^{\textsc{trip},{\lleno}}_{\textsc{ltc}}$, i.e. those trip arcs connecting nodes with the superscript $\lleno$.
Regarding the LTR formulation, PD$_1$ and PD$_2$ can be adapted only for instances with a single truck and PD$_3$ when $V' = V$, since the variables $X$ are no longer indexed by trucks. 

\vspace{-5pt}
\paragraph{Sequencing inequalities}
There exists a minimum time for a truck to deliver a given set of requests.
Although it is a naive estimation since the restrictions on drivers are ignored, it can be used as a lower bound on the route duration.
For example, the truck needs to (i) travel from its start location to the pickup location of the first request, (ii) wait until the pickup time window opens, (iii) load the cargo, (iv) travel to the delivery location, (v) wait until the delivery time window opens, (vi) unload the cargo, (vi) travel to the pickup location of the second request, and so on until delivering the last request.
Observe that all possible sequences of requests must be considered since their order is not known in advance.
Despite the computation of this parameter involves solving a Pickup-and-Delivery Problem with Time Windows \citep[see][]{DUMAS19917}, few requests are tractable even by brute-force algorithms.

Given $R' \subset R$ with $|R'| \geq 2$ and $v \in V$, let $\text{dur}^{\textsc{d}}_0(R',v)$ be the minimum time needed by $v$ to deliver all requests in $R'$.
By noting that $v$ can only deliver $|R'|-1$ requests of $R$ in a time less than $\text{dur}^{\textsc{p}}_0(R',v)$, the following inequalities are obtained.
\begin{align}
\sum_{r \in R'} \hspace{10pt} \sum_{\mathclap{\substack{e \in E^{\textsc{d},r}_{\textsc{lt}}:\\ \fin(e) < \text{dur}^{\textsc{d}}_0(R',v)}}} X_{ve} & \leq |R'| - 1 & \forall\ R' \subseteq R,\ v \in V. \nonumber \tag{SEC$_1$}
\end{align}
Again, to avoid making the description too long, similar inequalities are included for the pickups (i.e. considering the minimum time needed by $v$ until picking up the last request of $R'$). 
Notice that this approach spans from the beginning of the planning horizon.
Alternative inequalities can be proposed when considering an arbitrary $i \in \mathcal{I}$ as the beginning.
In this case, let $\text{dur}^{\textsc{d}}(R',i)$ be the minimum time a truck needs to deliver all requests in $R'$ beginning at $i$.
It is important to remark that, as the location of a truck at a particular instant is unknown in advance, it does not count the initial trip to the pickup location of the first request, and as a consequence, is independent of the truck.
\begin{align}
\sum_{r \in R'} \hspace{10pt} \sum_{\mathclap{\substack{e \in E^{\textsc{d},r}_{\textsc{lt}}:\\ \inicio(e) \geq i,\\ \fin(e) < i + \text{dur}^{\textsc{d}}(R',i)}}} X_{ve} & \leq |R'| - 1 & \forall\ R' \subseteq R,\ v \in V,\ i \in \mathcal{I}. \nonumber \tag{SEC$_2$}
\end{align}
Again, they can be adapted to the LTR formulation only for instances with a single truck.

\section{Computational experiments}\label{sec.computational.experiments}

In this section, the performance of the previous ILP formulations is compared and the incorporation of the proposed valid inequalities as cutting planes is evaluated, using a B\&C algorithm from a commercial ILP solver.
Experiments are performed by a computer equipped with Ubuntu 18.04 64bits, 6 GB of memory, and a core i7-9700 at 3.00 GHz.
It is common to use the benchmark program \emph{fmax} and the benchmark instance \emph{r500.5} (both available in \citet{DIMACS}) to allow future comparisons of our results with different machines. 
In our case, 2.86 s (user time) were spent.
The algorithms are implemented in \texttt{C++11} and call CPLEX Optimization Studio 20.1.0 through the Concert Technology API \citep{CPLEX}.
The CPLEX parameters have the default values, except a time limit of 2 h per instance, a single thread, and a deterministic search mode.

\begin{figure}
    \centering
    \includegraphics[scale=0.33]{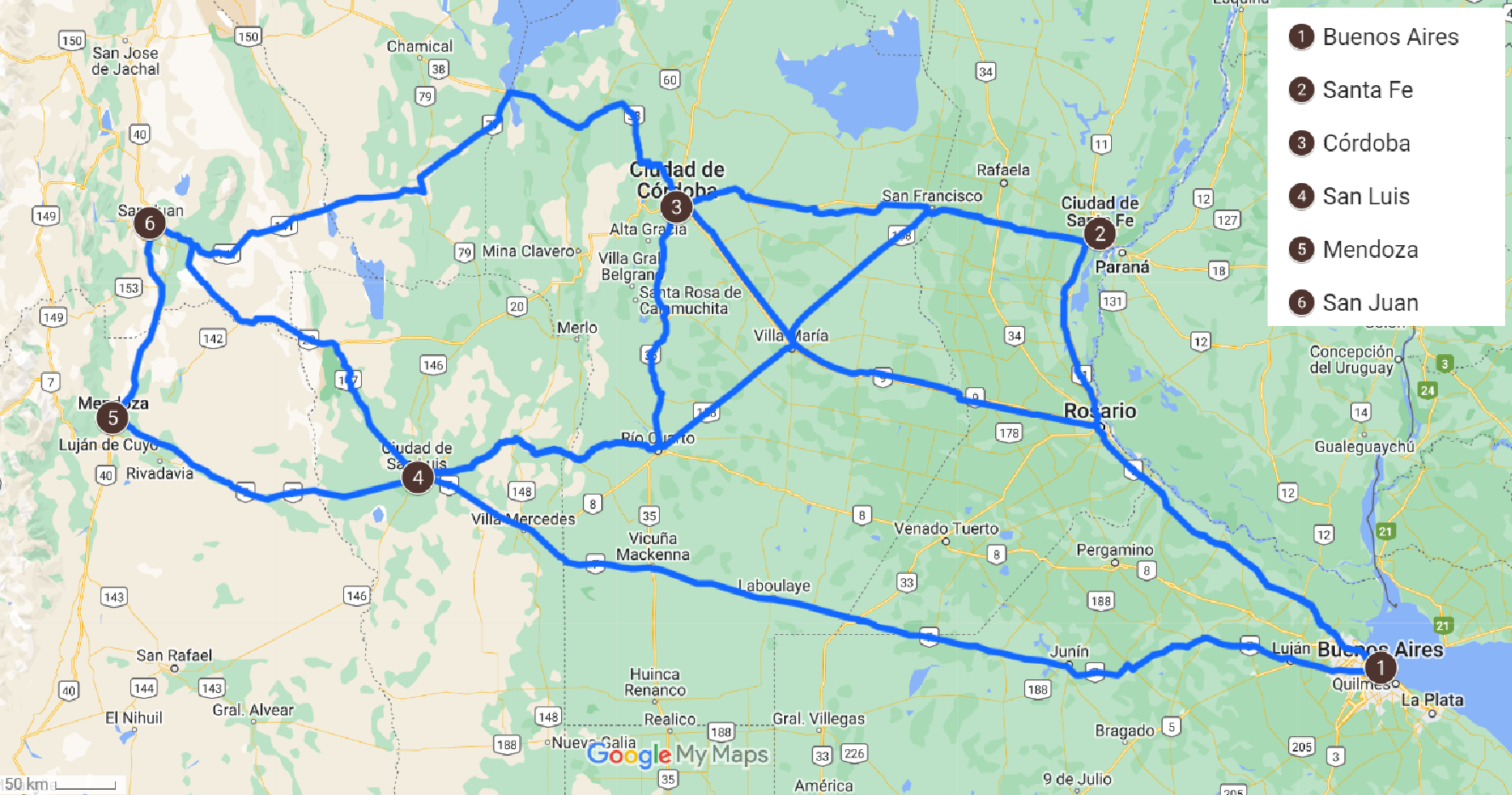}
    \caption{Fragment of the Argentine road network; extracted from Google Maps.}
    \label{fig.road.network}
\end{figure}

For this purpose, the following sets of random instances are generated, based on the road network depicted in Fig. \ref{fig.road.network}.
The generation considers a 1-hour discretisation of the planning horizon, a speed of 90 km/h for trucks and taxis, a travel cost for trucks equal to the travel time, and a travel cost for taxis twice that of trucks (since taxis must return to their origin).
\begin{itemize}
\item \texttt{S1}: 15 instances with 3 locations ($L = \{1,2,3\}$), $|P| \in \{4,5,6\}$, $H = 7$, $|V| = 1$ and $|D| = 2$. 
\item \texttt{S2}: 15 instances with $|P| \in \{7,8,9\}$ and the other parameters as \texttt{S1}. 
\item \texttt{S3}: 15 instances with 3 locations, $|P| \in \{8,9,10\}$, $|V| = 2$, and $|D| = 4$.
\item \texttt{S4}: 15 instances with 6 locations and the other parameters as \texttt{S1}.
\item \texttt{S5}: 4 instances with 3 locations, $|P| \in \{40,42\}$, $H \in \{40,42\}$, $|V| = 1$, and $|D| = 3$. 
\end{itemize}

\begin{figure}
\begin{subfigure}{0.49\hsize}
\includegraphics[width=\hsize]{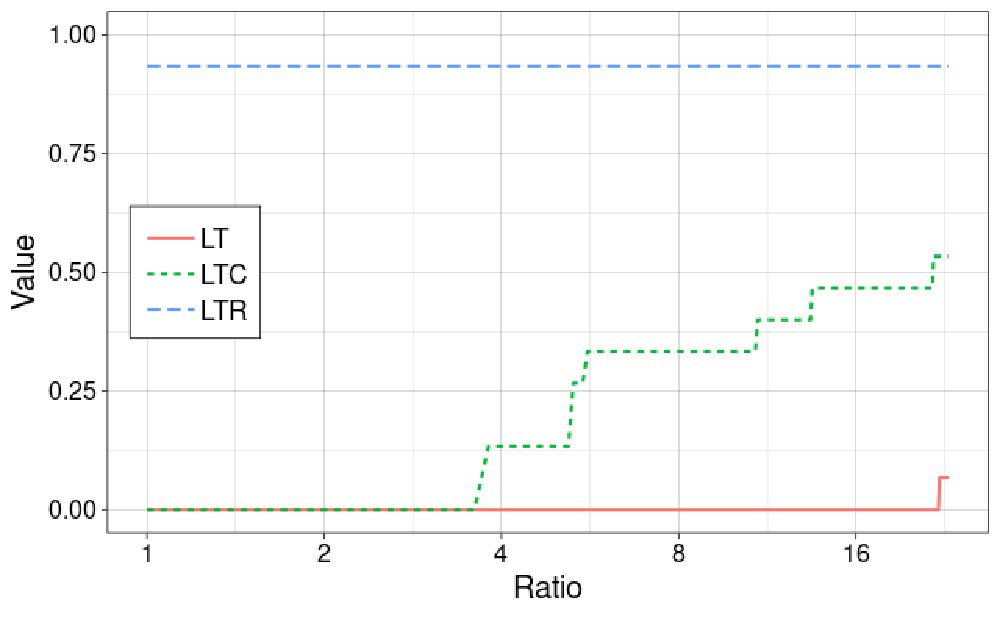}
\caption{Performance profiles of the algorithms.}
\label{fig.ilp.formulations.performance.profile}
\end{subfigure}
\hfill
\begin{subfigure}{0.49\hsize}
\includegraphics[width=\hsize]{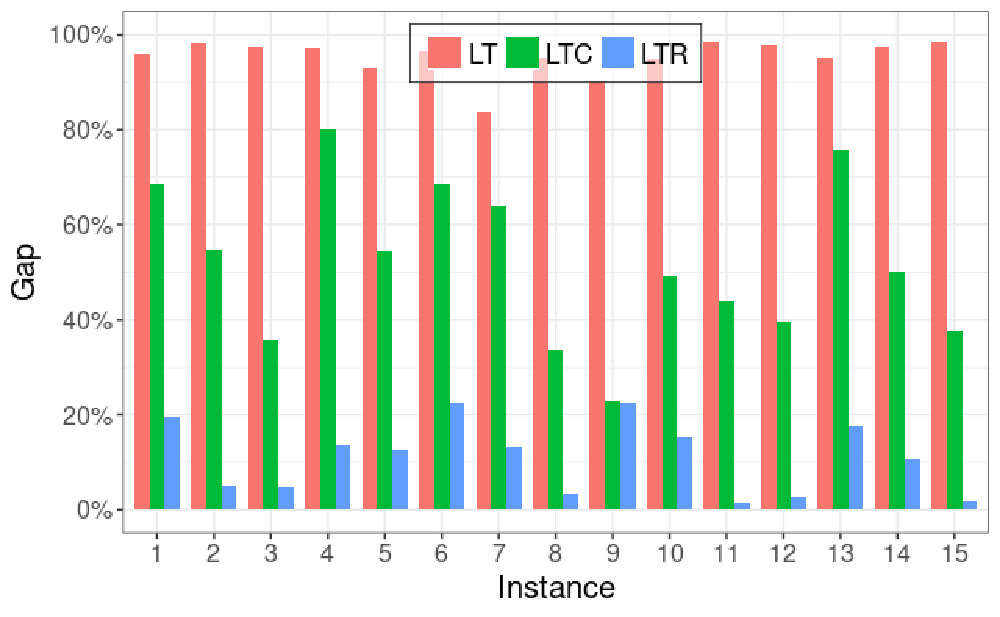}
\caption{Gap between the optimal value of the formulations and their LRs.}
\label{fig.ilp.formulations.gap}
\end{subfigure}
\caption{Comparison of ILP formulations when solving instances of \texttt{S1} (part II).}
\label{fig.ilp.formulations}
\end{figure}

\vspace{-5pt}
\paragraph{First experiment} For the first experiment, the ILP formulations, namely LT, LTC, and LTR, presented in Section \ref{sec.ilp}, are directly tackled with CPLEX.
The results are shown in Table \ref{tab.ilp.formulations}, where the first columns report the name of the formulation, the average number of variables and constraints, and the total number of solved instances (to optimality).
After that, the average execution time in seconds over all instances (7200 s is considered for unsolved instances) and over the solved instances are displayed.
Then, the average relative gap over all instances (0 \% is considered for solved instances) and over the unsolved instances are reported.
The last column shows the average number of explored nodes over the solved instances. 
The performance profile of the algorithms is depicted in Fig. \ref{fig.ilp.formulations.performance.profile}, using the execution time as the performance measure \citep[see][]{PerfProfile}.
Fig. \ref{fig.ilp.formulations.gap} shows, for each instance, the relative gap between the optimal value of the ILP formulation and its linear relaxation (LR).
It is worth mentioning that the optimal value is known for all the instances of \texttt{S1}, as will be shown in the second experiment. 
The average execution time of the LRs is 2.17 s (LT), 2.88 s (LTC), and 5.17 s (LTR).

\begin{table}
\centering
\caption{Comparison of ILP formulations when solving instances of \texttt{S1} (part I).}
\label{tab.ilp.formulations}
\begin{tabular*}{\textwidth}{@{}@{\extracolsep{\fill}}lllllllll@{}}
\hline
Form. & Var. & Cons. & Solved & \multicolumn{2}{l}{Time (s)} & \multicolumn{2}{l}{Gap (\%)} & Nodes \\[-6pt]
& & & & \multicolumn{2}{l@{}}{\hrulefill} & \multicolumn{2}{l@{}}{\hrulefill} & \multicolumn{1}{l@{}}{\hrulefill} \\
& & & & All & Sol. & All & Unsol. & Sol. \\
\hline
LT & 8682 & 5813 & 1 & 7109 & 5828 & 49 & 53 & 139372 \\
LTC & 10166 & 6151 & 8 & 4597 & 2319 & 14 & 33 & 9154 \\
LTR & 16102 & 7828 & 14 & 894 & 444 & 0.8 & 11 & 1915 \\
\hline
\end{tabular*}
\end{table}

A big difference is observed in the behaviour of the algorithms, mainly in the number of instances solved, where LTR far outperforms the others and LTC far outperforms LT.
The results show that the formulations are more competitive as the digraphs have greater structure, despite the increase in the number of variables and constraints.
At the same time, the LRs are tighter as the digraphs have greater structure, which might explain the difference observed.
Consequently, the next experiments focus on improving the best two algorithms (based on LTC and LTR) to solve larger instances.

\begin{figure}
\begin{subfigure}{0.49\hsize}
\includegraphics[width=\hsize]{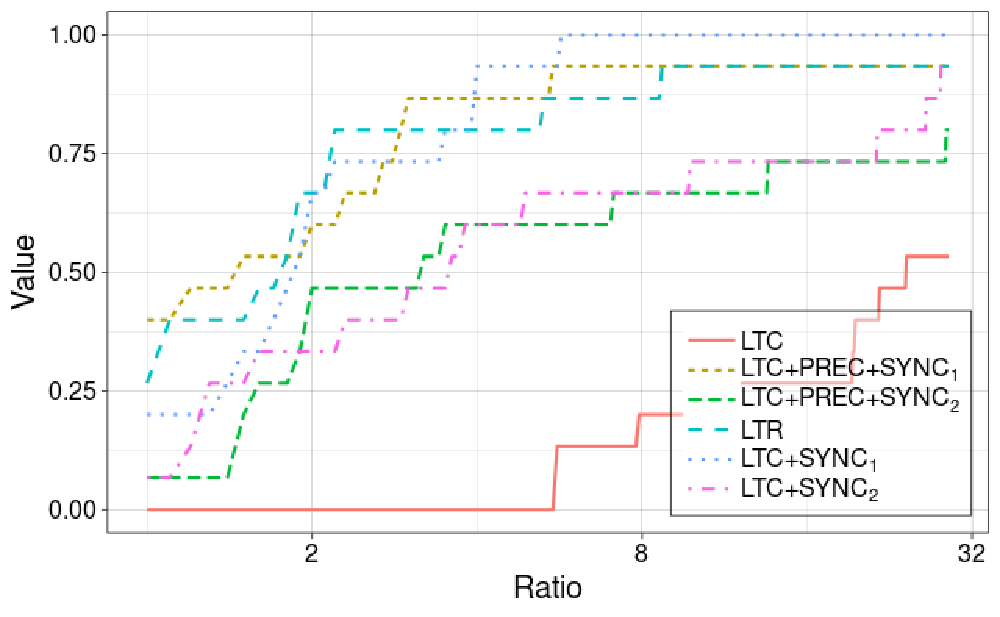}
\caption{Performance profiles of the algorithms.}
\label{fig.alt.ilp.formulations.performance.profile}
\end{subfigure}
\hfill
\begin{subfigure}{0.49\hsize}
\includegraphics[width=\hsize]{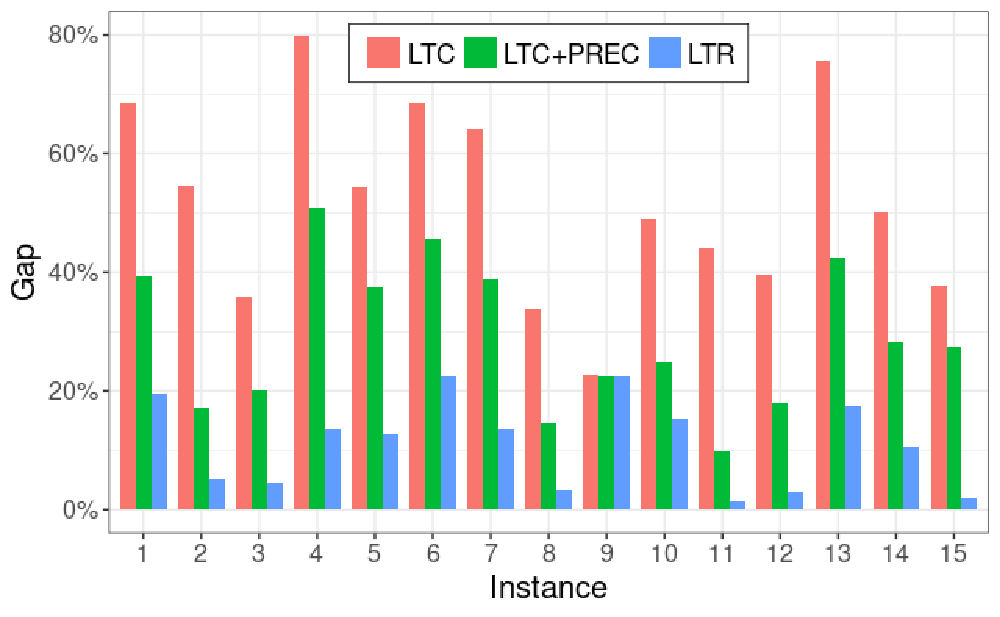}
\caption{Gap between the optimal value of the formulations and their LRs.}
\label{fig.alt.ilp.formulations.gap}
\end{subfigure}
\caption{Comparison of alternative ILP formulations when solving instances of \texttt{S1} (part II).}
\label{fig.alt.ilp.formulations}
\end{figure}

\vspace{-5pt}
\paragraph{Second experiment}
The next experiment evaluates the alternative formulations for LTC and LTR mentioned at the beginning of Section \ref{sec.valid.ineq}.
Their names are augmented with the family of constraints being replaced, e.g. in ``LTC + PREC + SYNC$_1$'', constraints (\ref{form.LT.LTX.Truck.Precedence}) and those in (\ref{form.LT.LTX.Synchronisation}) for pickup and delivery arcs are replaced with PREC and SYNC$_1$, respectively.
It is worth remembering that the former family is already implied in LTR, so LTR + PREC does not make sense in the current experiment.
The results are shown in Table \ref{tab.alt.ilp.formulations}, where the average number of variables can be seen from the previous table.
Fig. \ref{fig.alt.ilp.formulations.performance.profile} presents the performance profiles, where some of the less performant formulations are not shown to facilitate reading.
The LRs are analysed in Fig. \ref{fig.alt.ilp.formulations.gap}, where constraints SYNC$_1$ are SYNC$_2$ are omitted because they preserve the gaps.
The average execution time of the LRs is 2.57 s for LTC + PREC.

\begin{table}
\centering
\caption{Comparison of alternative ILP formulations when solving instances of \texttt{S1} (part I).}
\label{tab.alt.ilp.formulations}
\begin{tabular*}{\textwidth}{@{}@{\extracolsep{\fill}}llllllll@{}}
\hline
Form. & Cons. & Solved & \multicolumn{2}{l}{Time (s)} & \multicolumn{2}{l}{Gap (\%)} & Nodes \\[-6pt]
& & & \multicolumn{2}{l@{}}{\hrulefill} & \multicolumn{2}{l@{}}{\hrulefill} & \multicolumn{1}{l@{}}{\hrulefill} \\
& & & All & Sol. & All & Unsol. & Sol. \\
\hline
LTC  & 6151 & 8 & 4597 & 2319 & 14 & 33 & 9154\\
LTC + PREC & 6151 & 14 & 1274 & 851 & 1.5 & 22 & 7151\\
LTC + SYNC$_1$ & 5404 & 8 & 4674 & 2463 & 12 & 31 & 105495\\
LTC + SYNC$_2$ & 5404 & 8 & 4597 & 2319 & 7.3 & 16 & 168612\\
LTC + PREC + SYNC$_1$ & 5404 & 14 & 1019 & 577 & 2.6 & 39 & 4611\\
LTC + PREC + SYNC$_2$ & 5404 & 12 & 2201 & 951 & 0.6 & 3.1 & 54125\\ \hline
LTR & 7828 & 14 & 894 & 444 & 0.8 & 11 & 1915\\
LTR + SYNC$_1$ & 7082 & 15 & 917 & 917 & 0 & -- & 2318\\
LTR + SYNC$_2$ & 7082 & 14 & 2057 & 1690 & 1.3 & 19 & 56299\\
\hline
\end{tabular*}
\end{table}

A remarkable improvement is noticed, mainly in the number of solved instances, when LTC uses PREC, which greatly increases its competitiveness against LTR.
Again, the tighter LRs might explain this behaviour.
The best algorithms are based on LTC + PREC + SYNC$_1$, LTR, and LTR + SYNC$_1$.
Focusing on solving the greatest number of instances, we decided to continue working with LTC + PREC + SYNC$_1$ and LTR + SYNC$_1$, which from now on will directly be called LTC and LTR, respectively.

\vspace{-5pt}
\paragraph{Third experiment}
The following experiment analyses the impact of performing a warm start, i.e. providing the ILP solver with an external initial solution.
In our implementation, we use the metaheuristics developed in \citet{Lucci2023} to rapidly generate a heuristic solution, which is translated into a feasible solution for the integer linear programs. 
For each instance of \texttt{S1}, the average execution time of this initialisation is 7 s, and for the other sets it requires 18 s (\texttt{S2} and \texttt{S3}), 8 s (\texttt{S4}), and 464 s (\texttt{S5}).

\begin{figure}
\begin{subfigure}{0.49\hsize}
\includegraphics[width=\hsize]{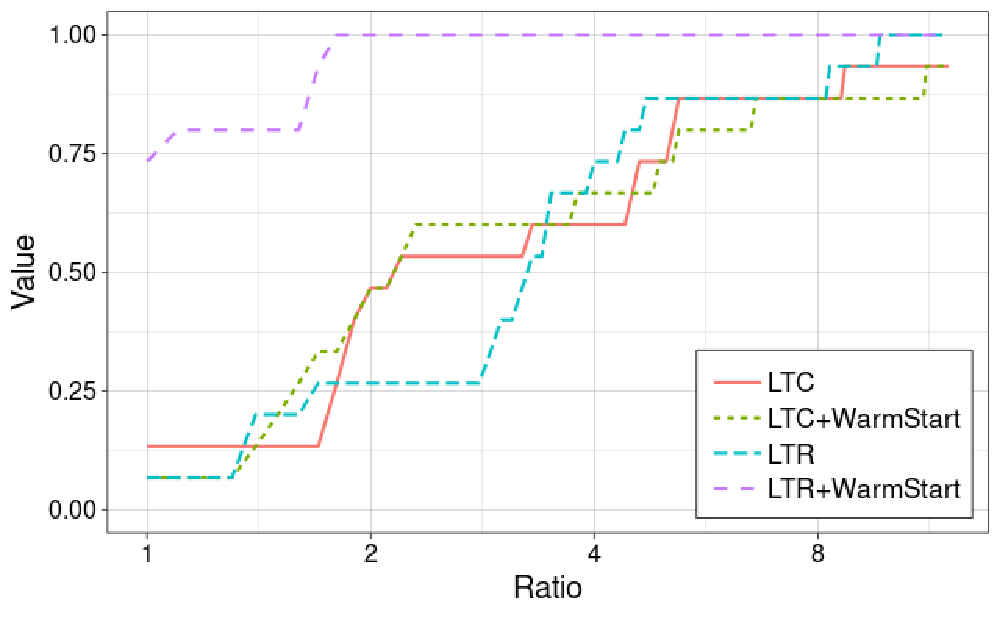}
\caption{Performance profiles of the algorithms.}
\label{fig.init.heur.performance.profile}
\end{subfigure}
\hfill
\begin{subfigure}{0.49\hsize}
\includegraphics[width=\hsize]{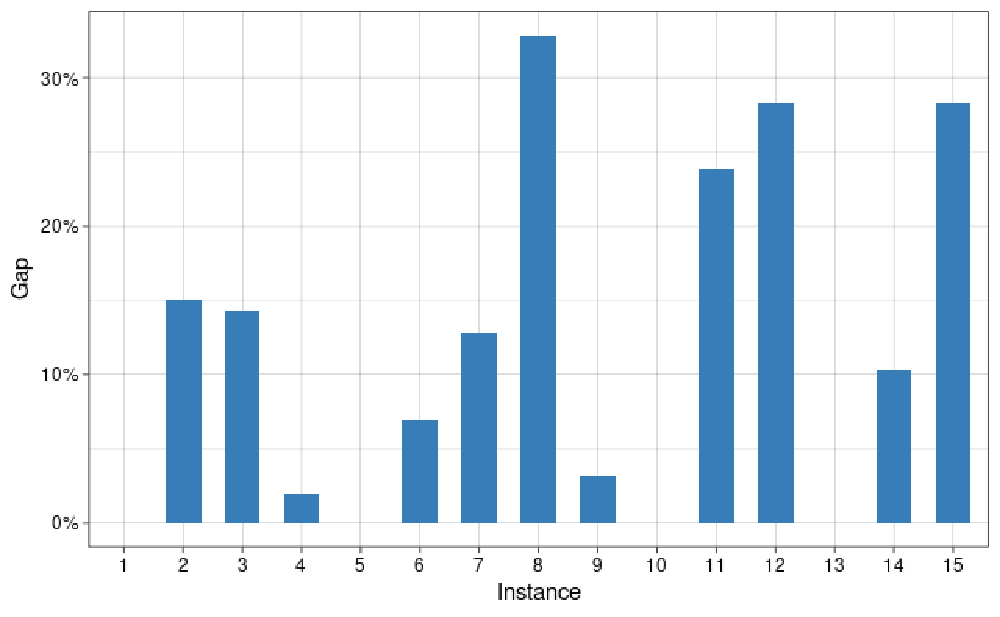}
\caption{Gap between the optimal value and the initial value.}
\label{fig.init.heur.gap}
\end{subfigure}
\caption{Impact of a warm start when solving instances of \texttt{S1}.}
\label{fig.init.heur}
\end{figure}

The results obtained for \texttt{S1} are discussed below. 
The warm start does not modify the number of solved instances but impacts the execution time.
In particular, LTR is notably faster when performing a warm start, whose average execution time over the solved instances decreases from  917 s to 647 s.
This is less noticeable for LTC, from 577 s to 484 s, but the gap for the only unsolved instance decreases from 39\% to 12\%.
To accompany the analysis, Fig. \ref{fig.init.heur.performance.profile} depicts the performance profiles, and Fig. \ref{fig.init.heur.gap} shows, for each instance, the gap between the optimal value and the value of the initial solution provided.
Since the warm start has proved to be useful, it will continue to be enabled in the remainder of this work.

\vspace{-5pt}
\paragraph{Fourth experiment}
It remains to evaluate the impact of the valid inequalities presented in Section \ref{sec.valid.ineq} for the tightening of the LRs.
In this experiment, we focus on the root node of the B\&C tree, while the remaining nodes are left for later.
The proposed methodology considers each family separately and incorporates all the valid inequalities into the models before solving their LRs.

It is worth noting that PD$_3$, SEC$_1$, and SEC$_2$ grow exponentially with the size of the input.
For this reason, a maximum size is introduced for the subsets of trucks ($V'$) and requests ($R'$).
The names of the inequalities are augmented accordingly, e.g., for $k \in \mathbb{N}$, ``PD$_3$-$V_k$'' considers all the inequalities of PD$_3$ such that $|V'| \leq k$ and similarly ``SEC$_1$-$R_k$'' demands $|R'| \leq k$.
However, this restriction might not be sufficient, and additional subfamilies are considered for PD$_3$ by limiting the initial time for loading and unloading.
In ``PD$_3$-$V_k$-A'', the start time of $e_1$ must coincide with the opening of a pickup time window and that of $e_2$ with the closing of a delivery time window.
For example, for time windows as Fig. \ref{fig.example.time.window}a, $e_1$ needs to start at 8:00 (on day 0 or 1) and $e_2$ at 20:00 (on day 0 or 1).
On the other hand, ``PD$_3$-$V_k$-B'' is similar but with the conjunction of both conditions, e.g. $e_1$ at 8:00 or $e_2$ at 20:00.

\begin{figure}
\begin{subfigure}{0.49\hsize}
\includegraphics[width=\hsize]{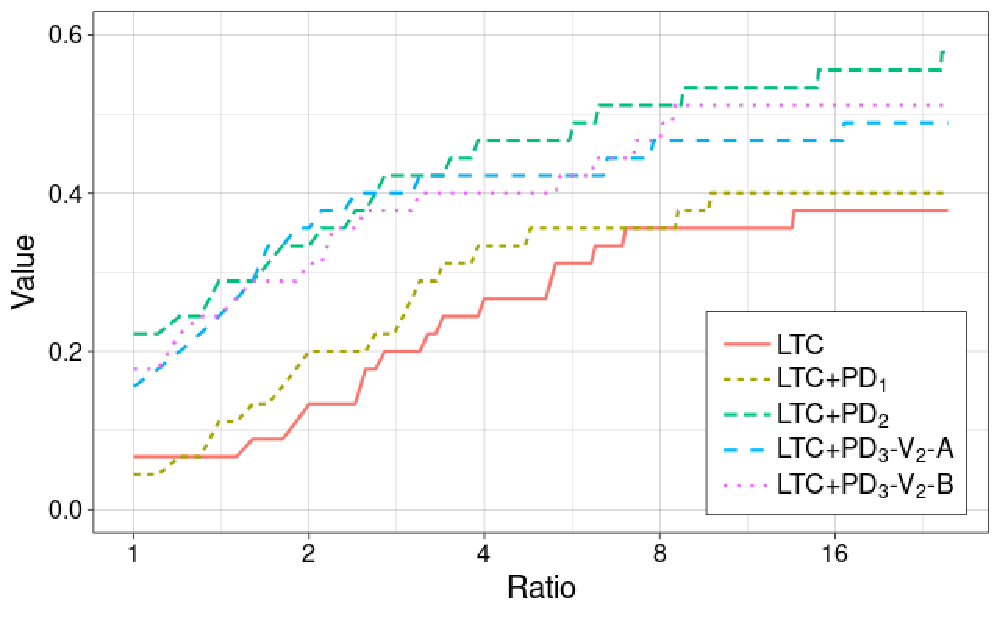}
\caption{Performance profiles of the algorithms.}
\label{fig.valid.ineq.pd.gap.performance.profile}
\end{subfigure}
\hfill
\begin{subfigure}{0.49\hsize}
\includegraphics[width=\hsize]{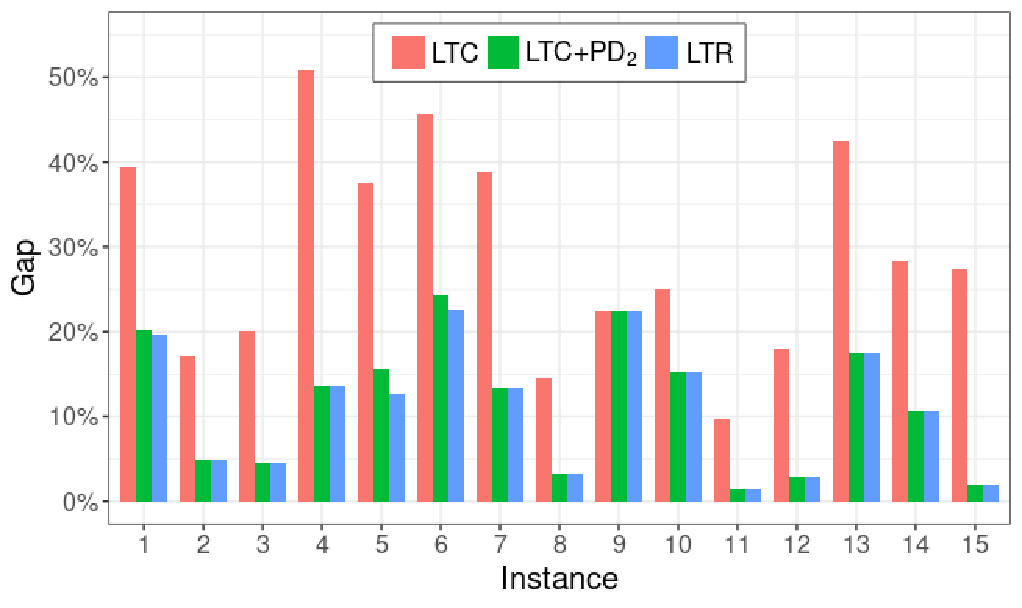}
\caption{Gap between the optimal value of the formulations and their LRs.}
\label{fig.valid.ineq.pd.gap}
\end{subfigure}
\caption{Impact of  PD inequalities on the entire tree for instances of \texttt{S2}, \texttt{S3}, and \texttt{S4} (part II).}
\label{fig.valid.ineq.pd}
\end{figure}

Tables \ref{tab.valid.ineq.root.I} and \ref{tab.valid.ineq.root.II} show the results for the instances of \texttt{S2} and \texttt{S3}, which are larger than \texttt{S1}.
The columns report the average number of valid inequalities, the average execution time, and the average relative gap with respect to the value of the best incumbent solution for the integer programs (which will be discovered later).
Families that have reduced the gap little or nothing (less than 0.2 percentage points) are excluded, e.g. PREC and PD for LTR in Table \ref{tab.valid.ineq.root.I}, SEC for LTC in Table \ref{tab.valid.ineq.root.II}.
Families that cannot be adapted are also omitted, e.g. most of the families for LTR in Table \ref{tab.valid.ineq.root.II}.
It is evident that the PD inequalities are highly effective for LTC.
In particular, the gap achieved by LTC + PD$_2$ is only 0.5 percentage points worse than LTR for the instances of \texttt{S2}.
Fig. \ref{fig.valid.ineq.pd.gap} exhibits a more detailed comparison of the gaps obtained for each instance of \texttt{S1}.
When the instances have more than one truck, the best value for $k$ in PD$_3$-$V_k$ seems to be $k = |V|$, but they are still outperformed by PD$_2$.
Regarding the SEC inequalities,
larger values of $k$ obtain better gaps, but the decrease is less abrupt from $k > 4$, while the number of inequalities and the execution times increase in counterpart.

\begin{table}
\centering
\caption{Impact of valid inequalities on the root node for instances of \texttt{S2} (part I).}
\label{tab.valid.ineq.root.I}
\begin{tabular*}{\textwidth}{@{\extracolsep{\fill}}llllllllll}
& \multicolumn{4}{l}{\hrulefill} & \multicolumn{4}{l}{\hrulefill} \\
& Form. & \#Ineq. & Time (s) & Gap (\%) & Algo. & \#Ineq. & Time (s) & Gap (\%)\\[-3pt]
& \multicolumn{4}{l}{\hrulefill} & \multicolumn{4}{l}{\hrulefill} \\
& LTC & 0 & 4.0 & 30.4 & LTR & 0 & 8.1 & 6.7\\
& LTC + PD$_1$ & 6 & 7.4 & 8.0 & LTR + SEC$_1$-$R_2$ & 30 & 8.7 & 6.7\\
& LTC + PD$_2$ & 686 & 9.5 & 7.2 & LTR + SEC$_1$-$R_3$ & 107 & 9.4 & 6.1\\
& LTC + PD$_3$-$V_1$-A & 349 & 7.1 & 9.8 & LTR + SEC$_1$-$R_4$ & 223 & 8.2 & 5.5\\
& LTC + PD$_3$-$V_1$-B & 8531 & 54.3 & 8.2 & LTR + SEC$_1$-$R_5$ & 335 & 8.3 & 5.4\\
& LTC + SEC$_1$-$R_2$ & 30 & 3.9 & 30.2 & LTR + SEC$_1$-$R_6$ & 406 & 8.5 & 5.4\\
& LTC + SEC$_1$-$R_3$ & 107 & 4.1 & 29.4 & LTR + SEC$_2$-$R_2$ & 6934 & 13.7 & 6.2\\
& LTC + SEC$_1$-$R_4$ & 223 & 4.5 & 28.9 & LTR + SEC$_2$-$R_3$ & 23233 & 30.7 & 6.0\\
& LTC + SEC$_1$-$R_5$ & 335 & 4.6 & 28.7 & LTR + SEC$_2$-$R_4$ & 45670 & 33.4 & 5.2\\
& LTC + SEC$_1$-$R_6$ & 406 & 4.9 & 28.7 & LTR + SEC$_2$-$R_5$ & 65948 & 36.9 & 5.1\\
& LTC + SEC$_2$-$R_2$ & 6934 & 6.5 & 29.8 & LTR + SEC$_2$-$R_6$ & 65948 & 36.9 & 5.1\\
& LTC + SEC$_2$-$R_3$ & 23233 & 11.8 & 29.3 & \multicolumn{4}{l}{\hrulefill} \\ % LTR + TAXI$_1$ & 3 & 8.1 & 6.1\\
& LTC + SEC$_2$-$R_4$ & 45670 & 14.1 & 28.4 & & & & \\ % LTR + TAXI$_2$ & 236 & 9.0 & 6.1\\
& LTC + SEC$_2$-$R_5$ & 65948 & 19.0 & 28.2 & & & &  \\
& LTC + SEC$_2$-$R_6$ & 78140 & 22.6 & 28.1 &  &  &  & \\[-3pt]
%& LTC + TAXI$_1$ & 3 & 3.9 & 29.6 &  &  &  & \\
%& LTC + TAXI$_2$ & 236 & 4.4 & 29.6 &  &  &  & \\[-3pt]
& \multicolumn{4}{l}{\hrulefill} & & & \\
\end{tabular*}
\end{table}

\begin{table}
\centering
\caption{Impact of valid inequalities on the root node for instances of \texttt{S3} (part II).}
\label{tab.valid.ineq.root.II}
\begin{tabular*}{0.5\textwidth}{@{\extracolsep{\fill}}llll}
\hline
Form. & Ineq. & Time (s) & Gap (\%)\\
\hline
LTC & 0 & 7.4 & 39\\
LTC + PD$_1$ & 12 & 14.6 & 24\\
LTC + PD$_2$ & 1403 & 24.9 & 5.6\\
LTC + PD$_3$-$V_1$-A & 789 & 10.3 & 38\\
LTC + PD$_3$-$V_1$-B & 18273 & 142.3 & 38\\
LTC + PD$_3$-$V_2$-A & 1184 & 17.4 & 15\\
LTC + PD$_3$-$V_2$-B & 27409 & 147.7 & 12\\
\hline
\end{tabular*}
\end{table}

Next, some of the best-performing families at the root node are selected and evaluated on the entire tree.
Small families are managed as an initial part of the constraint set of the models, thus all the valid inequalities are present in each node.
Otherwise, larger families are grouped into a pool of cuts and added to each node by the ILP solver on demand.

The results for PD inequalities are reported in Table \ref{tab.valid.ineq.pd} and Fig. \ref{fig.valid.ineq.pd.gap.performance.profile}.
In addition, the column ``Mode'' tells how the valid inequalities are handled, i.e. as part of the initial constraint set (``Init.'') or through a pool of cuts (``Pool''), and ``Cuts'' is the average number of cuts added from the pool per instance.
It is clear that PD$_2$, PD$_3$-$V_2$-$A$, and PD$_3$-$V_2$-$B$ conduct to a higher number of solved instances, lower execution times, and tighter gaps; in particular, the former slightly outperforms the last two.

\begin{table}
\centering
\caption{Impact of PD inequalities on the entire tree for instances of \texttt{S2}, \texttt{S3}, and \texttt{S4} (part I).}
\label{tab.valid.ineq.pd}
\begin{tabular*}{\textwidth}{@{\extracolsep{\fill}}lllllllllll}
\hline
Algo. & Mode & Ineq. & Cons. & Solved & \multicolumn{2}{l}{Time (s)} & \multicolumn{2}{l}{Gap (\%)} & Nodes & Cuts \\[-6pt]
& & & & & \multicolumn{2}{l@{}}{\hrulefill} & \multicolumn{2}{l@{}}{\hrulefill} & \multicolumn{1}{l@{}}{\hrulefill} & \\
& & & & & All & Sol. & All & Unsol. & Sol. & \\
\hline
LTC & -- & 0 & 10100 & 17 & 5474 & 2632 & 18 & 29 & 8415 & --\\
LTC + PD$_1$ & Init. & 10 & 10110 & 18 & 5330 & 2525 & 13 & 21 & 42325 & --\\
LTC + PD$_2$ & Init. & 912 & 11013 & 26 & 4459 & 2456 & 6.3 & 15 & 25916 & --\\
LTC + PD$_3$-$V_2$-A & Init. & 587 & 10687 & 22 & 4632 & 1946 & 9.5 & 19 & 6206 & --\\
LTC + PD$_3$-$V_2$-B & Pool & 13837 & 10100 & 23 & 4614 & 2140 & 8.5 & 17 & 5065 & 379\\
\hline
\end{tabular*}
\end{table}

\begin{figure}
\begin{subfigure}{0.49\hsize}
\includegraphics[width=\hsize]{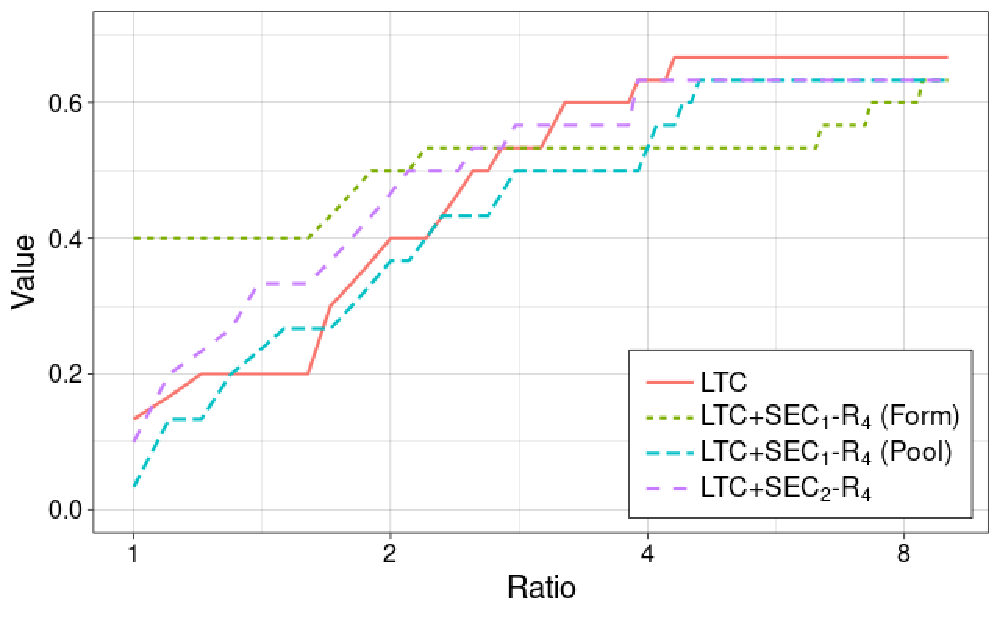}
\caption{Performance profiles of the algorithms based on LTC.}
\label{fig.valid.ineq.SEC.LTC}
\end{subfigure}
\hfill
\begin{subfigure}{0.49\hsize}
\includegraphics[width=\hsize]{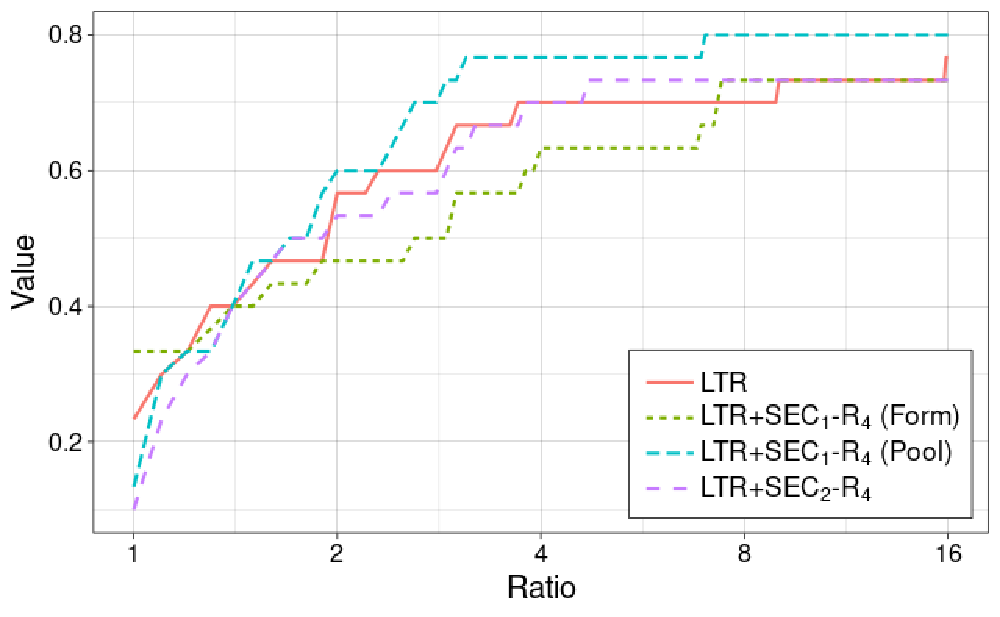}
\caption{Performance profiles of the algorithms based on LTR.}
\label{fig.valid.ineq.SEC.LTR}
\end{subfigure}
\caption{Impact of SEC inequalities on the entire tree for instances of \texttt{S2} and \texttt{S4} (part II).}
\label{fig.valid.ineq.SEC}
\end{figure}

While keeping PD$_2$ as part of the constraint set of LTC, SEC inequalities are now evaluated on the entire tree.
The results are reported in Table \ref{tab.valid.ineq.SEC} for the instances of \texttt{S2} and \texttt{S4} (\texttt{S3} is not considered since these inequalities were not effective on the root node when there is more than one truck).
The performance profiles of the algorithms based on LTC and LTR can be seen in Fig. \ref{fig.valid.ineq.SEC.LTC} and \ref{fig.valid.ineq.SEC.LTR}, respectively.
The results are inconclusive for LTC as no algorithm seems to outperform the others, e.g. the number of solved instances decreases.
Instead for LTR, SEC$_1$-$R_4$ (Pool) seems superior as the number of solved instances increases and the performance profile has the most wins within a ratio of 1.5 or higher.
It is interesting to mention that, although SEC$_1$-$R_4$ (Init) has the most wins within a ratio of 1, the performance decays notoriously as the ratio increases.
In addition, very few cuts are added from the pool.

\begin{table}
\centering
\caption{Impact of SEC inequalities on the entire tree for instances of \texttt{S2} and \texttt{S4} (part I).}
\label{tab.valid.ineq.SEC}
\small
\begin{tabular*}{\textwidth}{@{\extracolsep{\fill}}llllllllllll}
\hline
Algo. & Mode & Ineq. & Var. & Cons. & Solved & \multicolumn{2}{l}{Time (s)} & \multicolumn{2}{l}{Gap (\%)} & Nodes & Cuts \\[-6pt]
& & & & & & \multicolumn{2}{l@{}}{\hrulefill} & \multicolumn{2}{l@{}}{\hrulefill} & \multicolumn{1}{l@{}}{\hrulefill} & \\
& & & & & & All & Sol. & All & Unsol. & Sol. & \\
\hline
LTC & -- & 0 & 18906 & 10788 & 20 & 3989 & 2384 & 7.4 & 22 & 28398 & -- \\
LTC + SEC$_1$-$R_4$ & Init. & 129 & 18906 & 10917 & 19 & 3882 & 1961 & 7.7 & 21 & 13447 & -- \\
LTC + SEC$_1$-$R_4$ & Pool & 129 & 18906 & 10788 & 19 & 4201 & 2464 & 7.9 & 22 & 19197 & 0.5 \\
LTC + SEC$_2$-$R_4$ & Pool & 26427 & 18906 & 10788 & 19 & 4049 & 2225 & 8.1 & 22 & 13601 & 3.1 \\
\hline
LTR & -- & 0 & 32664 & 13477 & 23 & 2652 & 1268 & 3.6 & 15 & 4834 & -- \\
LTR + SEC$_1$-$R_4$ & Init. & 129 & 32664 & 13606 & 22 & 2963 & 1422 & 2.9 & 11 & 42406 & -- \\
LTR + SEC$_1$-$R_4$ & Pool & 129 & 32664 & 13477 & 24 & 2452 & 1265 & 2.4 & 12 & 4008 & 0.5 \\
LTR + SEC$_2$-$R_4$ & Pool & 26427 & 32664 & 13477 & 22 & 2548 & 857 & 3.2 & 12 & 3537 & 2.3 \\
\hline
\end{tabular*}
\end{table}

\vspace{-5pt}
\paragraph{Fifth experiment}
As a final experiment, the best algorithms for LTC and LTR are compared against each other, considering the instances of \texttt{S2}, \texttt{S3}, and \texttt{S4}.
To summarise the results, LTR shows a clear superiority in the number of solved instances, 37 vs. 26.
Every instance solved by LTC is also solved by LTR and more quickly (except two of them), with an average execution time of 734s vs. 2250s.
Furthermore, when none of them prove optimality, LTR reports better gaps (except for one instance), with an average gap of 12\% vs. 19\%.
In 9 of the 19 unsolved instances by LTC, the final lower bound is at least 5\% worse than LTR, and in 3 of them (all from \texttt{S4}), at least 24\%.
%To accompany this analysis, Fig. \ref{fig.final.comparison} depicts the performance profiles.
%
%\begin{figure}
%\centering
%\includegraphics[width=0.49\hsize]{exp5.eps}
%\caption{Final comparison of the algorithms based on LTC and LTR for instances of \texttt{S2}, \texttt{S3}, and \texttt{S4}.}
%\label{fig.final.comparison}
%\end{figure}

Then, the limits of the algorithms are analised using the instances of \texttt{S5}, which have the greatest number of requests and the longest planning horizon.
In this opportunity, we migrate to a more powerful computer, using the 16 cores of an i7-10700 at 2.9GHz (parallel optimization), 6 GB of memory, and a time limit of 10 hours.
On the one hand, LTC has an average of 158k variables and 95k constraints.
The LRs are successfully solved on the root node but no other nodes are explored within the time limit.
The LRs are tightened with 3 to 5 rounds of general-purpose cutting planes and the final gaps vary from 6.7\% to 18.5\%.
On the other hand, LTR runs out of memory while writing the models for 3 instances and during the presolve for the other, with 485k variables and 186k constraints.

Finally, we reconsider the first algorithm based on LTC, without alternative constraints or valid inequalities, but keeping the warm start (with the same initial solutions) for the instances of \texttt{S5}.
Again, only the LRs are solved on the root node and there is not enough time to branch. 
Although more rounds of cutting planes are performed, 7 to 11, the final gaps deteriorate greatly, varying from 69.8\% y 79.4\%.

%% Conclusion of the paper
\section{Conclusions}\label{sec.conclusions}

In this work, we addressed a SVRCSP in long-haul transport.
The goal was to minimise the travel costs and the delay penalties involved in transporting pickup-and-delivery requests with time windows over a multi-day planning horizon.
Many real-world requirements were considered, such as the hour-of-service regulation of the drivers. 
Unlike most approaches in the literature, the correspondence between trucks and drivers was not fixed and could be exchanged in some locations.
%In addition, drivers are allowed to travel as passengers in the company's trucks or use an external taxi service for an additional cost.
Different digraphs were defined to represent truck and driver routes separately as directed paths, where nodes were indexed by location and time.
Some of them incorporated greater structure (additional nodes and arcs) to ensure that certain forbidden paths could not exist by construction.
Three compact ILP formulations were proposed to model the problem based on these digraphs, namely LT, LTC, and LTR.
These models allowed truck and driver routes to be easily synchronised in time and space.
Among the most notable differences, the number of variables of LT and LTC strongly depended on the number of trucks, whereas LTR on the number of requests.
Many families of valid inequalities were presented and some of them were exponential in size.

Extensive computational experiments were performed on random instances to compare the ILP formulations and evaluate specific valid inequalities as cutting planes.
One of the most interesting conclusions was that, as the digraphs gained structure and the formulations contained fewer families of constraints, the linear relaxations were tighter and the number of solved instances increased.
Furthermore, tightening LTC with precedence inequalities and pickup-and-delivery trip inequalities produced notable performance improvements, primarily in the number of solved instances.
In contrast, the algorithm based on LTR, which initially performed best, supported notable reductions in execution time with the addition of a warm start and a pool of cuts with the sequencing inequalities.
A final comparison revealed that the latter significantly outperformed the former in instances with 3-6 locations, 7-10 requests, 1-2 trucks, 2-4 drivers, and a 1-week planning horizon.
However, the algorithm based on LTC stood out for providing relative gaps for larger instances, with 40-42 requests and a planning horizon of 40-42 days.

In conclusion, this study has contributed to the development of exact methods for solving SVRCSPs.
To our knowledge, this is the first time that solutions with proven optimality have been obtained for the case study.
We believe that these methods can be extended to other similar studies, but potential adaptations might be necessary to handle different hour-of-service regulations.
This research also suggested some topics that would be interesting to study in the future.
Relaxing some of the requirements might achieve further cost improvements, e.g. allowing drivers to share taxis.
Another line of research is the definition of new valid inequalities, which might help to reduce the differences in the performance of the algorithms, considering the limitation of LTR when distinguishing trucks.
Regarding the implementation part, it would be interesting to develop specific separation routines for some families of valid inequalities, as well as primal heuristics.
Finally, since the models generally have many more variables than constraints, a branch-and-price scheme could be competitive to address larger instances.

% Acknowledgements
\section*{Acknowledgments}
This work was partially supported by grants PIP-1900 (CONICET) and PICT-2020-03032 (ANPCyT).

% References

\nocite{*}

\bibliographystyle{apalike}
\bibliography{lucci-arxiv}

% Appendix
% \vspace*{-10pt}
% \appsection{Appendix A}\label{App:A}

% \subsection{Appendix subhead}

\end{document}